\documentclass[leqno,12pt]{article}  
\usepackage{amsmath}
\usepackage{amssymb}

\usepackage[german,english]{babel}
\usepackage{amsthm}

\title{On Families of Pure Slope $L$-Functions}

\author{\textsc{Elmar Grosse-Kl\"onne}}
\date{}

\theoremstyle{plain} 
\newtheorem{satz}{Theorem}[section]  
\newtheorem{kor}[satz]{Corollary}  
\newtheorem{lem}[satz]{Lemma}  
\newtheorem{pro}[satz]{Proposition}  
\newcommand{\ho}{\mbox{\rm Hom}}  
 
\newcommand{\spec}{\mbox{\rm Spec}}  

\newcommand{\spm}{\mbox{\rm Sp}}  

\newcommand{\ord}{\mbox{\rm ord}}  
\newcommand{\ke}{\mbox{\rm Ker}}  
\newcommand{\koke}{\mbox{\rm Coker}}  

\newcommand{\tr}{\mbox{\rm Tr}}

\newcommand{\sym}{\mbox{\rm Sym}}
\newcommand{\id}{\mbox{\rm id}}

\theoremstyle{remark}

\theoremstyle{definition}

\begin{document}
\maketitle
\footnote[0]{}
\footnote[0]{\textit{Key words and phrases}. Unit root L-function, Dwork's conjecture, $\sigma$-module, slope decomposition, weight space.}

\footnote[0]{  I wish to express my sincere thanks to Robert Coleman and Daqing Wan. Manifestly this work heavily builds on ideas of them, above all on Wan's limiting module construction. Wan invited me to begin further elaborating his methods, and directed my attention to many interesting problems involved. Coleman asked me for the meromorphic continuation to the whole character space and provided me with some helpful notes \cite{cormk}. In particular the important functoriality result \ref{funclim} for the limiting module and the suggestion of varying it rigid analytically is due to him. Thanks also to Matthias Strauch for discussions on the weight space.                          }

\begin{abstract}

Let $R$ be the ring of integers in a finite extension $K$ of $\mathbb{Q}_p$, let $k$ be its residue field and let $\chi:\pi_1(X)\to R^{\times}=GL_{1}(R)$ be a "geometric" rank one representation of the arithmetic fundamental group of a smooth affine $k$-scheme $X$. We show that the locally $K$-analytic characters $\kappa:R^{\times}\to\mathbb{C}_p^{\times}$ are the $\mathbb{C}_p$-valued points of a $K$-rigid space ${\cal W}$ and that $$L(\kappa\circ\chi,T)=\prod_{\overline{x}\in X}\frac{1}{1-(\kappa
\circ\chi)(Frob_{\overline{x}})T^{\deg(\overline{x})}},$$viewed as a two variable function in $T$ and $\kappa$, is meromorphic on $\mathbb{A}_{\mathbb{C}_p}^1\times{\cal W}$. On the way we prove, based on a construction of Wan, a slope decomposition for ordinary overconvergent (finite rank) $\sigma$-modules, in the Grothendieck group of nuclear $\sigma$-modules.

\end{abstract}

%


\begin{center} {\bf Introduction} 
\end{center}

In a series of remarkable papers \cite{waann} \cite{warko} \cite{wahrk}, Wan recently proved a long outstanding conjecture of Dwork on the $p$-adic meromorphic continuation of unit root $L$-functions arising from an ordinary family of algebraic varieties defined over a finite field $k$. We begin by illustrating his result by a concrete example. Fix $n\ge0$ and let $Y$ be the affine $n+1$-dimensional $\mathbb{F}_p$-variety in $\mathbb{A}^1\times\mathbb{G}_m^{n+1}$ defined by $$z^p-z=x_0+\ldots+\ldots x_n.$$Define $u:Y\to\mathbb{G}_m$ by sending $(z,x_0,\ldots,x_n)$ to $x_0x_1\cdots x_n$. For $r\ge1$ and $y\in\mathbb{F}_{p^r}^{\times}$ let $Y_y/\mathbb{F}_{p^r}$ be the fibre of $u$ above $y$. For $m\ge 1$ let $Y_y(\mathbb{F}_{p^{rm}})$ be the set of $\mathbb{F}_{p^{rm}}$-rational points and $(Y_y)_0$ the set of closed points of $Y_y/\mathbb{F}_{p^r}$ (a closed point $z$ is an orbit of an $\overline{\mathbb{F}}_{p^r}$-valued point under the $p^r$-th power Frobenius map $\sigma_{p^r}$; its degree $\deg_r(z)$ is the smallest positive integer $d$ such that $\sigma_{p^r}^d$ fixes the orbit pointwise). The zeta function of $Y_y/\mathbb{F}_{p^r}$ is$$Z(Y_y/\mathbb{F}_{p^r},T)=\exp(\sum_{m=1}^{\infty}\frac{|Y_y(\mathbb{F}_{p^{rm}})|}{m}T^m)=\prod_{z\in(Y_y)_0}\frac{1}{1-T^{\deg_r(z)}}.$$On the other hand for a character $\Psi:\mathbb{F}_p\to\mathbb{C}$ define the Kloosterman sum$$K_{m}(y)=\sum_{\stackrel{x_i\in\mathbb{F}_{p^{rm}}^{\times}}{x_0x_1\cdots x_n=y}}\Psi(\tr_{\mathbb{F}_{p^{rm}}/\mathbb{F}_{p}}(x_0+x_1+\ldots+x_n))$$and let $L_{\Psi}(Y,T)$ be the series such that$$T{\mbox {dlog}} L_{\Psi}(y,T)=\sum_{m=1}^{\infty}K_m(y)T^m.$$Then, as series,$$\prod_{\Psi}L_{\Psi}(Y,T)=Z(Y_y/\mathbb{F}_{p^r},T),$$hence to understand $Z(Y_y/\mathbb{F}_{p^r},T)$ we need to understand all the $L_{\Psi}(y,T)$. Suppose $\Psi$ is non-trivial. It is known that $L_{\Psi}(y,T)$ is a polynomial of degree $n+1$: there are algebraic integers $\alpha_0(y),\ldots,\alpha_n(y)$ such that$$L_{\Psi}(y,T)^{(-1)^{n-1}}=(1-\alpha_0(y)T)\cdots(1-\alpha_n(y)T).$$These $\alpha_i(y)$ have complex absolute value $p^{rn/2}$ and are $\ell$-adic units for any prime $\ell\ne p$. We ask for their $p$-adic valuation and their variation with $y$. Embedding $\overline{\mathbb{Q}}\to\overline{\mathbb{Q}_p}$ we have $\alpha_i(y)\in\mathbb{Q}_p(\pi)$ where $\pi^{p-1}=-p$. Sperber has shown that we may order the $\alpha_i(y)$ such that $\ord_p(\alpha_i(y))=i$ for any $0\le i\le n$. Fix such an $i$ and for $k\in\mathbb{Z}$ consider the $L$-function$$\prod_{y\in(\mathbb{G}_m)_0/\mathbb{F}_p}\frac{1}{1-(p^{-i}\alpha_i(y))^kT^{\deg_1(y)}}$$(here $\deg_1(y)$ is the minimal $r$ such that $y\in\mathbb{F}_{p^r}^{\times}$, and $(\mathbb{G}_m)_0/\mathbb{F}_p$ is the set of closed points of $\mathbb{G}_m/\mathbb{F}_p$ defined similarly as before). A priori this series defines a holomorphic function only on the open unit disk. Dwork conjectured and Wan proved that it actually extends to a meromorphic function on $\mathbb{A}_{\mathbb{C}_p}^1$, and varies uniformly with $k$ in some sense. Now let ${\cal W}$ be the rigid space of locally $\mathbb{Q}_p(\pi)$-analytic characters $\kappa$ of the group of units in the ring of integers of $\mathbb{Q}_p(\pi)$. In this paper we show that $$L(T,\kappa)=\prod_{y\in(\mathbb{G}_m)_0/\mathbb{F}_p}\frac{1}{1-\kappa(p^{-i}\alpha_i(y))T^{\deg_1(y)}}$$defines a meromorphic function on $\mathbb{A}_{\mathbb{C}_p}^1\times{\cal W}$. Specializing $\kappa\in{\cal W}$ to the character $r\mapsto r^k$ for $k\in \mathbb{Z}$ we recover Wan's result. The conceptual way to think of this example is in terms of $\sigma$-modules: $\mathbb{F}_p$ acts on $Y$ via $z\mapsto z+a$ for $a\in\mathbb{F}_p$. This induces an action of $\mathbb{F}_p$ on the relative $n$-th rigid cohomology $\mathbb{R}^nu_{rig,*}{\mathcal O}_Y$ of $u$, and over $\mathbb{Q}_p(\pi)$ the latter splits up into its eigencomponents for the various characters of $\mathbb{F}_p$. The $\Psi$-eigencomponent $(\mathbb{R}^nu_{rig,*}{\mathcal O}_Y)^{\Psi}$ is an overconvergent $\sigma$-module, and $L_{\Psi}(y,T)^{(-1)^{n-1}}$ is the characteristic polynomial of Frobenius acting on its fibre in $y$. Crucial is the slope decomposition of $(\mathbb{R}^nu_{rig,*}{\mathcal O}_Y)^{\Psi}$: it means that for fixed $i$ the $\alpha_i(y)$ vary rigid analytically with $y$ in some sense. We are thus led to consider Dwork's conjecture, i.e. Wan's theorem, in the following general context.

Let $R$ be the ring of integers in a finite extension $K$ of $\mathbb{Q}_p$, let $\pi$ be a uniformizer and $k$ the residue field. Let $X$ be a smooth affine $k$-scheme, let $A$ be the coordinate ring of a lifting of $X$ to a smooth affine weak formal $R$-scheme (so $A$ is a wcfg-algebra) and let $\widehat{A}$ be the $p$-adic completion of $A$. Let $\sigma$ be an $R$-algebra endomorphism of $A$ lifting the $q$-th power Frobenius endomorphism of $X$, where $q=|k|$. A finite rank $\sigma$-module over $\widehat{A}$ (resp. over $A$) is a finite rank free $\widehat{A}$-module (resp. ${A}$-module) together with a $\sigma$-linear endomorphism $\phi$. A finite rank $\sigma$-module over $\widehat{A}$ is called overconvergent if it arises by base change $A\to\widehat{A}$ from a finite rank $\sigma$-module over $A$. Let the finite rank overconvergent $\sigma$-module $\Phi$ over $\widehat{A}$ be ordinary, in the strong sense that it admits a Frobenius stable filtration such that on the $j$-th graded piece we have: the Frobenius is divisible by $\pi^j$ and multiplied with $\pi^{-j}$ it defines a unit root $\sigma$-module ${\Phi}_j$, i.e. a $\sigma$-module whose linearization is bijective. (Recall that unit root $\sigma$-modules over $\hat{A}$ are the same as continuous representations of $\pi_1(X)$ on finite rank free $R$-modules.) Although $\Phi$ is overconvergent, ${\Phi}_j$ will in general not be overconvergent; and this is what prevented Dwork from proving what is now Wan's theorem: the $L$-function $L({\Phi}_j,T)$ is meromorphic on ${\mathbb A}^1_{\mathbb{C}_p}$. Moreover he proved the same for powers (=iterates of the $\sigma$-linear endomorphism) ${\Phi}_j^k$ of ${\Phi}_j$ and showed that in case ${\Phi}_j$ is of rank one the family $\{L({\Phi}_j^k,T)\}_{k\in\mathbb{Z}}$ varies uniformly with $k\in\mathbb{Z}$ in a certain sense. At the heart of Wan's striking method lies his "limiting $\sigma$-module" construction which allows him to reduce the analysis of the not necessarily overconvergent ${\Phi}_j$ to that of overconvergent $\sigma$-modules --- at the cost of now working with overconvergent $\sigma$-modules of infinite rank, but which are nuclear. To the latter a generalization of the Monsky trace formula can be applied which expresses $L({\Phi}_j^k,T)$ as an alternating sum of Fredholm determinants of completely continuous Dwork operators.

The first aim of this paper is to further explore the significance of the limiting $\sigma$-module construction which we think to be relevant for the search of good $p$-adic coefficients on varieties in characteristic $p$. Following an argument of Coleman \cite{cormk} we give a functoriality result for this construction. This is then used to prove (Theorem \ref{hndecov}) a slope decomposition for ordinary overconvergent finite rank $\sigma$-modules, in the Grothendieck group $\Delta(\widehat{A})$ of nuclear $\sigma$-modules over $\widehat{A}$. More precisely, we show that any ${\Phi}_j$ as above, not necessarily overconvergent, can be written, in $\Delta(\widehat{A})$, as a sum of virtual nuclear {\it overconvergent} $\sigma$-modules. (This is the global version of the decomposition of the corresponding $L$-function found by Wan.) Our second aim is to strengthen Wan's uniform results on the family $\{L({\Phi}_j^k,T)\}_{k\in\mathbb{Z}}$ in case ${\Phi}_j$ is of rank one. More generally we replace ${\Phi}_j$ by the rank one unit root $\sigma$-module $\det({\Phi}_j)$ if ${\Phi}_j$ has rank $>1$. Let $\det{\Phi}_j$ be given by the action of $\alpha\in\widehat{A}^{\times}$ on a basis element. For $\overline{x}\in X$ a closed point of degree $f$ let $x:\widehat{A}\to R_f$ be its Teichm\"uller lift, where $R_f$ denotes the unramified extension of $R$ of degree $f$. Then $$ \alpha_{\overline{x}} = x(\alpha\sigma(\alpha)\ldots\sigma^{f-1}(\alpha)) $$ lies in $R^{\times}$. We prove that for any locally $K$-analytic character $\kappa:R^{\times}\to\mathbb{C}_p^{\times}$ the twisted $L$-function$$L(\alpha,T,\kappa)=\prod_{\overline{x}\in X}\frac{1}{1-\kappa(\alpha_{\overline{x}})T^{\deg(\overline{x})}}$$is $p$-adic meromorphic on $\mathbb{A}_{\mathbb{C}_p}^1$, and varies rigid analytically with $\kappa$. More precisely, building on work of Schneider and Teitelbaum \cite{scte}, we use Lubin-Tate theory to construct a smooth $\mathbb{C}_p$-rigid analytic variety ${\cal W}$ whose $\mathbb{C}_p$-valued points are in natural bijection with the set $\ho_{K\mbox{-}an}(R^{\times},\mathbb{C}_p^{\times})$ of locally $K$-analytic characters of $R^{\times}$. Then our main theorem is:

\begin{satz}\label{hauptsatz} On the $\mathbb{C}_p$-rigid space $\mathbb{A}_{\mathbb{C}_p}^1\times{\cal W}$ there exists a meromorphic function $L_\alpha$ whose pullback to $\mathbb{A}_{\mathbb{C}_p}^1$ via $\mathbb{A}_{\mathbb{C}_p}^1\to\mathbb{A}_{\mathbb{C}_p}^1\times{\cal W},\, t\mapsto(t,\kappa)$, for any $\kappa\in\ho_{K\mbox{-}an}(R^{\times},\mathbb{C}_p^{\times})={\cal W}(\mathbb{C}_p)$ is a continuation of $L(\alpha,T,\kappa)$.\\
\end{satz}

The statement as formulated in the abstract above follows by the well known correspondence between representations of the fundamental group and unit-root $\sigma$-modules.

The analytic variation of the $L$-series $L(\alpha,T,\kappa)$ with the weight $\kappa$ makes it meaningful to vastly generalize the eigencurve theme studied by Coleman and Mazur \cite{eigenc} in connection with the Gouv\^{e}a-Mazur conjecture. Namely, we can ask for the divisor of the two variable meromorphic function $L_{\alpha}$ on $\mathbb{A}_{\mathbb{C}_p}^1\times{\cal W}$. From a general principle in \cite{cobmf} we already get: for fixed $\lambda\in\mathbb{R}_{>0}$, the difference between the numbers of poles and zeros of $L_{\alpha}$ on the annulus $|T|=\lambda$ is locally constant on ${\cal W}$. We hope for better qualitative results if the $\sigma$-module over $A$ giving rise to the $\sigma$-module $\Phi$ over $\widehat{A}$ carries an overconvergent integrable connection, i.e. is an overconvergent $F$-isocrystal on $X$ in the sense of Berthelot. The eigencurve from \cite{eigenc} comes about in this context as follows: The Fredholm determinant of the $U_p$-operator acting on overconvergent $p$-adic modular forms is a product of certain power rank one unit root $L$-functions arising from the universal ordinary elliptic curve, see \cite{cobmf}. Also, again in the general case, the $p$-adic $L$-function on ${\cal W}$ which we get by specializing $T=1$ in $L_{\alpha}$ should be of particular interest.

The proof of Theorem \ref{hauptsatz} consists of two steps. First we prove (this is essentially Corollary \ref{wafa}) the meromorphic continuation to ${\mathbb A}_{\mathbb{C}_p}^{1}\times{\cal W}^0$ for a certain open subspace ${\cal W}^0$ of ${\cal W}$ which meets every component of ${\cal W}$: the subspace of characters of the type $\kappa(r)=r^{\ell}u(r)^x$ for $\ell\in\mathbb{Z}$ and small $x\in\mathbb{C}_p$, with $u(r)$ denoting the one-unit part of $r\in R^{\times}$. (In particular, ${\cal W}^0$ contains the characters $r\mapsto \kappa_k(r)=r^k$ for $k\in\mathbb{Z}$; for these we have $L({\Phi}_j^k,T)=L(\alpha,T,\kappa_k)$.) For this we include $\det({\Phi}_j)$ in a {\it family} of nuclear $\sigma$-modules, parametrized by ${\cal W}^0$: namely, the factorization into torsion part and one-unit part and then exponentiation with $\ell\in\mathbb{Z}$ resp. with small $x\in\mathbb{C}_p$ makes sense not just for $R^{\times}$-elements but also for $\alpha$, hence an analytic family of rank one unit root $\sigma$-modules parametrized by ${\cal W}^0$. In the Grothendieck group of ${\cal W}^0$-parametrized {\it families} of nuclear $\sigma$-modules, we write this deformation family of $\det({\Phi}_j)$ as a sum of virtual families of nuclear {\it overconvergent} $\sigma$-modules. In each fibre $\kappa\in{\cal W}^0$ we thus obtain, by an infinite rank version of the Monsky trace formula, an expression of the $L$-function $L(\alpha,T,\kappa)$ as an alternating product of characteristic series of nuclear Dwork operators. While this is essentially an "analytic family version" of Wan's proof (at least if $X=\mathbb{A}^n$), the second step, the extension to the whole space $\mathbb{A}_{\mathbb{C}_p}^1\times{\cal W}$, needs a new argument. We use a certain integrality property (w.r.t. ${\cal W}$) of the coefficients of (the logarithm of) $L_{\alpha}$ which we play out against the already known meromorphic continuation on ${\mathbb A}_{\mathbb{C}_p}^{1}\times{\cal W}^0$. However, we are not able to extend the limiting modules from ${\cal W}^0$ to all of ${\cal W}$; as a consequence, for $\kappa\in{\cal W}-{\cal W}^0$ we have no interpretation of $L(\alpha,T,\kappa)$ as an alternating product of characteristic series of Dwork operators. Note that for $K=\mathbb{Q}_p$, the locally $K$-analytic characters of $R^{\times}=\mathbb{Z}_p^{\times}$ are precisely the {\it continuous} ones; the space ${\cal W}^0$ in that case is the weight space considered in \cite{cobmf} while ${\cal W}$ is that of \cite{eigenc}. 

Now let us turn to some technical points. Wan develops his limiting $\sigma$-module construction and the Monsky trace formula for nuclear overconvergent infinite rank $\sigma$-modules only for the base scheme $X=\mathbb{A}^n$. General base schemes $X$ he embeds into $\mathbb{A}^n$ and treats (the pure graded pieces of) finite rank overconvergent $\sigma$-modules on $X$ by lifting them with the help of Dwork's $F$-crystal to $\sigma$-modules on $\mathbb{A}^n$ having the same $L$-functions. We work instead in the infinite rank setting on arbitrary $X$. Here we need to overcome certain technical difficulties in extending the finite rank Monsky trace formula to its infinite rank version. The characteristic series through which we want to express the $L$-function are those of certain Dwork operators $\psi$ on spaces of overconvergent functions with {\it non fixed} radius of overconvergence. To get a hand on these $\psi$'s one needs to write these overconvergent function spaces as direct limits of appropriate affinoid algebras on which the restrictions of the $\psi$'s are completely continuous. Then statements on the $\psi$'s can be made if these affinoid algebras have a common system of orthogonal bases. Only for $X=\mathbb{A}^n$ we find such bases; but we show how one can pass to the limit also for general $X$. An important justification for proving the trace formula in this form (on general $X$, with function spaces with {\it non fixed} radius of overconvergence) is that in the future it will allow us to make full use of the overconvergent connection in case the $\sigma$-module over $A$ giving rise to the $\sigma$-module $\Phi$ over $\widehat{A}$ underlies an overconvergent $F$-isocrystal on $X$ (see above) --- then the limiting module also carries an overconvergent connection. Deviating from \cite{waann} \cite{warko}, instead of working with formally free nuclear $\sigma$-modules with fixed formal bases we work, for concreteness, with the infinite square matrices describing them. This is of course only a matter of language.

A brief overview. In section 1 we show the existence of common orthogonal bases in overconvergent ideals which might be of some independent interest. In section 2 we define the $L$-functions and prove the trace formula. In section 3 we introduce the Grothendieck group of nuclear $\sigma$-modules (and their deformations). In section 4 we concentrate on the case where ${\phi}_j$ is the unit root part of $\phi$ and is of rank one: here we need the limiting module construction. In section 5 we introduce the weight space ${\cal W}$, in section 6 we prove (an infinite rank version of) Theorem \ref{hauptsatz}, and in section 7 (which logically could follow immediately after section 4) we give the overconvergent representation of ${\Phi}_j$.\\

{\it Notations:} By $|.|$ we denote an absolute value of $K$ and by $e\in\mathbb{N}$ the absolute ramification index of $K$. By $\mathbb{C}_p$ we denote the completion of a fixed algebraic closure of $K$ and by $\ord_{\pi}$ and $\ord_p$ the homomorphisms $\mathbb{C}_p^{\times}\to \mathbb{Q}$ with $\ord_{\pi}(\pi)=\ord_p(p)=1$. We write $\mathbb{N}_0=\mathbb{Z}_{\ge0}$. For $R$-modules $E$ with $\pi E\ne E$ we set$$\ord_{\pi}(x):=\sup\{r\in\mathbb{Q};\, r=\frac{n}{m} \mbox{ for some }n\in\mathbb{N}_0, m\in\mathbb{N}\mbox{ such that } x^m\in \pi^nE\}$$for $x\in E$. Similarly we define $\ord_p$ on such $E$. For $n\in \mathbb{N}$ we write $\mathbb{\mu}_n=\{x\in\mathbb{C}_p;\,x^n=1\}$. For an element $g$ in a free polynomial ring $A[X_1,\ldots,X_n]$ over a ring $A$ we denote by $\deg(g)$ its (total) degree. We will use the usual notations $|\alpha|=\sum_{i=1}^n\alpha_i$ for a multiindex $\alpha=(\alpha_0,\ldots,\alpha_n)\in\mathbb{N}_0^n$, and  $[r]\in\mathbb{Z}$ for a given $r\in\mathbb{Q}$: the unique integer with $[r]\le r<[r]+1$.\\

\section{Orthonormal bases of overconvergent ideals}
In this preparatory section we determine explicit orthonormal $K$-bases of ideals in overconvergent $K$-Tate algebras $T_n^c$ (\ref{onbas}). Furthermore we recall the complete continuity of certain Dwork operators (\ref{matlim}).\\
 
\addtocounter{satz}{1}{\bf \arabic{section}.\arabic{satz}} \newcounter{teenro1}\newcounter{teenro2}\setcounter{teenro1}{\value{section}}\setcounter{teenro2}{\value{satz}} For $c\in\mathbb{N}$ we let$$T_n^c:=\{\sum_{\alpha\in\mathbb{N}_0^n}b_{\alpha}\pi^{[\frac{|\alpha|}{c}]}X^{\alpha};\quad b_{\alpha}\in K, \lim_{|\alpha|\to\infty}|b_{\alpha}|=0\}.$$This is the ring of power series in $X_1,\ldots,X_n$ with coefficients in $K$, convergent on the polydisk  $$\{x\in\mathbb{C}_p^n;\quad\ord_{\pi}(x_i)\ge -\frac{1}{c}\mbox{ for all }1\le i\le n\}.$$ We view $T_n^c$ as a $K$-Banach module with the unique norm $|.|_c$ for which $\{\pi^{[\frac{|\alpha|}{c}]}X^{\alpha}\}_{\alpha\in\mathbb{N}_0^n}$ is an orthonormal basis (this norm is not power multiplicative). Suppose we are given elements $g_1,\ldots, g_r\in R[X_1,\ldots,X_n]-\pi R[X_1,\ldots,X_n]$. Let $\overline{g}_j\in k[X_1,\ldots,X_n]$ be the reduction of $g_j$, let $d_j=\deg(\overline{g}_j)\le\deg(g_j)$ be its degree.

\begin{lem}\label{normeins} For each $1\le j\le r$ and each $c>\max_j\deg(g_j)$ we have$$|\pi^{[\frac{|\alpha|+d_j}{c}]}X^{\alpha}g_j|_c=1.$$
\end{lem}

{\sc Proof:} Write $g_j=\sum_{\beta\in\mathbb{N}_0^n}b_{\beta}X^{\beta}$ with $b_{\beta}\in K$. There exists a $\beta_1\in\mathbb{N}_0^n$ with $|\beta_1|=d_j$ and $|b_{\beta_1}|=1$. Hence $$|\pi^{[\frac{|\alpha|+d_j}{c}]}X^{\alpha}b_{\beta_1}X^{\beta_1}|_c=|\pi^{[\frac{|\alpha|+|\beta_1|}{c}]}X^{\alpha+\beta_1}|_c=1.$$Now let $\beta\in\mathbb{N}_0^n$ be arbitrary, with $b_{\beta}\ne 0$. If $|\beta|>d_j$ then $|b_{\beta}|\le|\pi|$. Hence $$|\pi^{[\frac{|\alpha|+d_j}{c}]}X^{\alpha}b_{\beta}X^{\beta}|_c\le|\pi^{[\frac{|\alpha|+d_j}{c}]-[\frac{|\alpha|+|\beta|}{c}]+1}\pi^{[\frac{|\alpha|+|\beta|}{c}]}X^{\alpha+\beta}|_c.$$But $[\frac{|\alpha|+d_j}{c}]-[\frac{|\alpha|+|\beta|}{c}]+1\ge0$ because $b_{\beta}\ne0$, hence $c>|\beta|$. Thus, $$|\pi^{[\frac{|\alpha|+d_j}{c}]}X^{\alpha}b_{\beta}X^{\beta}|_c\le1.$$On the other hand, if $|\beta|\le d_j$, then $[\frac{|\alpha|+d_j}{c}]\ge[\frac{|\alpha|+|\beta|}{c}]$ and $|b_{\beta}|\le 1$, and again we find$$|\pi^{[\frac{|\alpha|+d_j}{c}]}X^{\alpha}b_{\beta}X^{\beta}|_c\le1.$$We are done.\\ 

\addtocounter{satz}{1}{\bf \arabic{section}.\arabic{satz}}\newcounter{defide1}\newcounter{defide2}\setcounter{defide1}{\value{section}}\setcounter{defide2}{\value{satz}} The Tate algebra in $n$ variables over $K$ is the algebra $$T_n:=\{\sum_{\alpha\in\mathbb{N}_0^n}b_{\alpha}X^{\alpha};\quad b_{\alpha}\in K, \lim_{|\alpha|\to\infty}|b_{\alpha}|=0\}.$$Let $I_{\infty}$ (resp. $I_c$) be the ideal in $T_n$ (resp. in $T_n^c$) generated by $g_1,\ldots, g_r$. As all ideals in $T_n^c$, the ideal $I_c$ is closed in $T_n^c$. We view $I_c$ as a $K$-Banach module with the norm $|.|_c$ induced from $T_n^c$.

\begin{lem}\label{exaide} If $I_{\infty}\subset T_n$ is a prime ideal, $I_{\infty}\ne T_n$, then $I_c=I_{\infty}\cap T_n^c$ for $c>>0$.
\end{lem} 

{\sc Proof:} For $c>>0$ also $I_c$ is a prime ideal in $T_n^c$. The open immersion of $K$-rigid spaces $\spm(T_n)\to\spm(T_n^c)$ induces an open immersion $V(I_{\infty})\to V(I_c)$ of the respective zero sets of $g_1,\ldots,g_r$. That $I_c$ is prime means that $V(I_c)$ is irreducible, and $I_{\infty}\ne T_n$ means that $V(I_{\infty})$ is non empty. Hence an element of $I_{\infty}\cap T_n^c$, since it vanishes on $V(I_{\infty})$, necessarily also vanishes on $V(I_c)$. By Hilbert's Nullstellensatz (\cite{bgr}) it is then an element of $I_c$.\\

Now we fix an integer $c'>\max_j\deg(g_j)$. By \ref{normeins} we find a subset $E$ of $\mathbb{N}_0^n\times\{1,\ldots,r\}$ such that $\{\pi^{[\frac{|\alpha|+d_j}{c'}]}X^{\alpha}g_j\}_{(\alpha,j)\in E}$ is an orthonormal basis of $I_{c'}$ over $K$.

\begin{satz}\label{onbas} For integers $c\ge c'$, the set $\{\pi^{[\frac{|\alpha|+d_j}{c}]}X^{\alpha}g_j\}_{(\alpha,j)\in E}$ is an orthonormal basis of $I_{c}$ over $K$.
\end{satz}

{\sc Proof:} Let $K^c$ be a finite extension of $K$ containing a $c$-th root $\pi^{\frac{1}{c}}$ and a $c'$-th root $\pi^{\frac{1}{c'}}$ of $\pi$. The absolute value $|.|$ extends to $K^c$. Any norm on a $K$-Banach module $M$ extends uniquely to a $K^c$-Banach module norm on $M\otimes_KK^c$, and we keep the same name for it. It is enough to show that $\{\pi^{[\frac{|\alpha|+d_j}{c}]}X^{\alpha}g_j\}_{(\alpha,j)\in E}$ is an orthonormal basis of $I_{c}\otimes_KK^c$ over $K^c$. Let $|.|'_c$ be the supremum norm on $T_n^c\otimes_KK^c$. This is the norm for which $\{X^{\frac{\alpha}{c}}\}_{\alpha\in\mathbb{N}_0^n}$ is an orthonormal basis over $K^c$. For $j\in\{1,\ldots,r\}$ write $g_j=\sum_{\beta\in\mathbb{N}_0^n}b_{\beta}X^{\beta}$ with $b_{\beta}\in K$. Then, by a computation similar to that in \ref{normeins} we find$$|\pi^{\frac{|\alpha|+d_j}{c}}X^{\alpha}b_{\beta}X^{\beta}|'_c=1\quad\mbox{ if }|\beta|=d_j\mbox{ and }|b_{\beta}|=1,$$$$|\pi^{\frac{|\alpha|+d_j}{c}}X^{\alpha}b_{\beta}X^{\beta}|'_c<1\quad\mbox{ otherwise.                 }\,\quad\quad\quad\quad\quad$$In particular it follows that $|\pi^{\frac{|\alpha|+d_j}{c}}X^{\alpha}g_j|'_c=1$. Now a comparison of expansions shows that $\{\pi^{[\frac{|\alpha|+d_j}{c}]}X^{\alpha}g_j\}_{(\alpha,j)\in E}$ is an orthonormal basis of $I_{c}\otimes_KK^c$ over $K^c$ with respect to $|.|_c$ if and only if $\{\pi^{\frac{|\alpha|+d_j}{c}}X^{\alpha}g_j\}_{(\alpha,j)\in E}$ is an orthonormal basis of $I_{c}\otimes_KK^c$ over $K^c$ with respect to $|.|'_c$. In particular it follows on the one hand that we only need to show that $\{\pi^{\frac{|\alpha|+d_j}{c}}X^{\alpha}g_j\}_{(\alpha,j)\in E}$ is an orthonormal basis of $I_{c}\otimes_KK^c$ over $K^c$ with respect to $|.|'_c$, and on the other hand it follows (applying the above with $c'$ instead of $c$) that $\{\pi^{\frac{|\alpha|+d_j}{c'}}X^{\alpha}g_j\}_{(\alpha,j)\in E}$ is an orthonormal basis of $I_{c'}\otimes_KK^c$ over $K^c$ with respect to $|.|'_{c'}$. Consider the isomorphism $$T_n^c\otimes_KK^c\cong T_n^{c'}\otimes_KK^c,\quad\quad\pi^{\frac{\alpha}{c}}X^{\alpha}\mapsto \pi^{\frac{\alpha}{c'}}X^{\alpha}$$which is isometric with respect to $|.|'_c$ resp. $|.|'_{c'}$. It does not necessarily map $I_c\otimes_KK^c$ to $I_{c'}\otimes_KK^c$. However, from our above computations of the values $|\pi^{\frac{|\alpha|+d_j}{c}}X^{\alpha}b_{\beta}X^{\beta}|'_c$ it follows that this isomorphism identifies the reductions of the elements of the set $\{\pi^{\frac{|\alpha|+d_j}{c}}X^{\alpha}g_j\}_{(\alpha,j)\in \mathbb{N}_0^n\times\{1,\ldots,r\}}$ with the reductions of the elements of the set $\{\pi^{\frac{|\alpha|+d_j}{c'}}X^{\alpha}g_j\}_{(\alpha,j)\in \mathbb{N}_0^n\times\{1,\ldots,r\}}$ (here by reduction we mean reduction modulo elements of absolute value $<1$). The $K^c$-vector subspaces spanned by these sets are dense in $I_c\otimes_KK^c$ resp. in $I_{c'}\otimes_KK^c$. Since for a subset of $|.|=1$ elements in an orthonormizable $K^c$-Banach module the property of being an orthonormal basis is equivalent to that of inducing an (algebraic) basis of the reduction, the theorem follows.\\  

\addtocounter{satz}{1}{\bf \arabic{section}.\arabic{satz}}\newcounter{ovcoc1}\newcounter{ovcoc2}\setcounter{ovcoc1}{\value{section}}\setcounter{ovcoc2}{\value{satz}} Let $B_K$ be a reduced $K$-affinoid algebra, i.e. a quotient of a Tate algebra $T_m$ over $K$ (for some $m$), endowed with its supremum norm $|.|_{\sup}$. Let $$B=(B_K)^0:=\{b\in B_K;\quad |b|_{\sup}\le 1\}.$$ For positive integers $m$ and $c$ let$$[m,c]:=[m,c]_{R}:=\{z\in R[[X_1,\ldots,X_n]];$$$$ z=\sum_{j=0}^{\infty}\pi^jp_j\mbox{ with }p_j\in R[X_1,\ldots,X_n]\mbox{ and }\deg(p_j)\le m+cj\}$$and $$[m,c]_{B}:=[m,c]\widehat{\otimes}_RB$$ (the $\pi$-adically completed tensor product). Note that for $m, c_1, c_2\in\mathbb{N}$ with $c_1<c_2$ we have $[m,c_1]_{B}\subset T_n^{c_2}\widehat{\otimes}_KB_K$ and also $(\cup_{m,c}[m,c]_{B})\otimes_RK=\cup_{c}(T_n^c\widehat{\otimes}_KB_K)$. Let$$R[X_1,\ldots,X_n]^{\dagger}:= R[X]^{\dagger}:=\bigcup_{m,c}[m,c].$$ Fix a Frobenius endomorphism $\sigma$ of $R[X]^{\dagger}$ lifting the $q$-th power Frobenius endomorphism of $k[X]$. Also fix a Dwork operator $\theta$ (with respect to $\sigma$) on $R[X]^{\dagger}$, i.e. an $R$-module endomorphism with $\theta(\sigma(x)y)=x\theta(y)$ for all $x,y\in R[X]^{\dagger}$. By \cite{fc3} 2.4 we have $\theta(T_n^c)\subset T_n^c$ for all $c>>0$, thus we get a $B_K$-linear endomorphism $\theta\otimes 1$ on $T_n^c\widehat{\otimes}_KB_K$.\\

\begin{pro}\label{matlim} Let $I$ be a countable set, $m',c'$ positive integers and ${\cal M}=(a_{i_1,i_2})_{i_1,i_2\in I}$ an $I\times I$-matrix with entries $a_{i_1,i_2}$ in $[m',c']_{B}$. Suppose that ${\cal M}$ is nuclear, i.e. that for each $M>0$ there are only finitely many $i_2\in I$ such that $\inf_{i_1}\ord_{\pi}a_{i_1,i_2}<M$. For $c>>0$ and $\beta\in\mathbb{N}_0^n$ develop $(\theta\otimes1)(\pi^{[\frac{|\beta|}{c}]}X^{\beta}a_{i_1,i_2})\in T_n^c\widehat{\otimes}_KB_K$ in the orthonormal basis $\{\pi^{[\frac{|\alpha|}{c}]}X^{\alpha}\}_{\alpha}$ of the $B_K$-Banach module $T_n^c\widehat{\otimes}_KB_K$ and let $G^{c}_{\{\alpha,i_1\}\{\beta,i_2\}}\in B_K$ for $\alpha\in\mathbb{N}_0^n$ be its coefficients:$$(\theta\otimes1)(\pi^{[\frac{|\beta|}{c}]}X^{\beta}a_{i_1,i_2})=\sum_{\alpha}G^{c}_{\{\alpha,i_1\}\{\beta,i_2\}}\pi^{[\frac{|\alpha|}{c}]}X^{\alpha}.$$Then for all $c>>0$ and all $M>0$ there are only finitely many pairs $(\alpha,i_2)\in \mathbb{N}_0^n\times I$ such that $$\inf_{\beta,i_1}\ord_{\pi}G^c_{\{\alpha,i_1\}\{\beta,i_2\}}<M.$$
\end{pro}

{\sc Proof:} For simplicity identify $I$ with $\mathbb{N}$. By \cite{fc3} 2.3 we find integers $r$ and $c_0$ such that $(\theta\otimes1)([qm,qc]_{B})\subset[m+r,c]_{B}$ for all $c\ge c_0$, all $m$. Increasing $c_0$ and $r$ we may assume that $a_{i_1,i_2}\in[q(r-1),c_0]_{B}$ for all $i_1,i_2$. Now let $c$ be so large that for $c'=c-1$ we have $qc'\ge c_0$. Then one easily checks that $X^{\beta}a_{i_1i_2}\in[q(r+[\frac{|\beta|}{q}]),qc']_{B}$ for all $\beta, i_1, i_2$. Hence $(\theta\otimes1)(X^{\beta}a_{i_1,i_2})\in[r+[\frac{|\beta|}{q}],c']$. This means $$|\alpha|\le r+[\frac{|\beta|}{q}]+c'(\ord_{\pi}(G^{c}_{\{\alpha,i_1\}\{\beta,i_2\}})+[\frac{|\alpha|}{c}]-[\frac{|\beta|}{c}])$$for all $\alpha$, and thus$$\ord_{\pi}(G^{c}_{\{\alpha,i_1\}\{\beta,i_2\}})\ge[\frac{|\beta|}{c}]-[\frac{|\alpha|}{c}]+\frac{|\alpha|-r-[\frac{|\beta|}{q}]}{c'}.$$Here the right hand side tends to infinity as $|\alpha|$ tends to infinity, uniformly for all $\beta$ --- independently of $i_1$ and $i_2$ --- because $c/q\le c'\le c$. Now let $M\in\mathbb{N}$ be given. By the above we find $N'(M)\in\mathbb{N}$ such that for all $\alpha$ with $|\alpha|\ge N'(M)$ we have $\ord_{\pi}(G^{c}_{\{\alpha,i_1\}\{\beta,i_2\}})\ge M$. Now fix $\alpha$. We have$$\ord_{\pi}(G^{c}_{\{\alpha,i_1\}\{\beta,i_2\}})\ge[\frac{|\beta|}{c}]-[\frac{|\alpha|}{c}]+\ord_{\pi}(\theta\otimes1)(X^{\beta}a_{i_1,i_2}).$$By nuclearity of ${\cal M}$ the right hand side tends to zero as $i_2$ tends to infinity, uniformly for all $i_1$, all $\beta$. In other words, there exists $N(\alpha,M)$ such that $\ord_{\pi}(G^{c}_{\{\alpha,i_1\}\{\beta,i_2\}})\ge M$ for all $i_2\ge N(\alpha,M)$, for all $i_1$, all $\beta$. Now set $$N(M)=N'(M)+\max\{N(\alpha,M);|\alpha|<N'(M)\}.$$Then we find $\inf_{\beta,i_1}\ord_{\pi}G^c_{\{\alpha,i_1\}\{\beta,i_2\}}\ge M$ whenever $|\alpha|+i_2\ge N(M)$. We are done.\\

\section{$L$-functions}

This section introduces our basic setting. We define nuclear (overconvergent) matrices (which give rise to nuclear (overconvergent) $\sigma$-modules), their associated $L$-functions and Dwork operators and give the Monsky trace formula (\ref{tracf}).\\ 

\addtocounter{satz}{1}{\bf \arabic{section}.\arabic{satz}}\newcounter{setgeo1}\newcounter{setgeo2}\setcounter{setgeo1}{\value{section}}\setcounter{setgeo2}{\value{satz}} Let $q\in\mathbb{N}$ be the number of elements of $k$, i.e. $k=\mathbb{F}_{q}$. Let $X=\spec(\overline{A})$ be a smooth affine connected $k$-scheme of dimension $d$. So $\overline{A}$ is a smooth $k$-algebra. By \cite{elk} it can be represented as $\overline{A}=A/\pi A$ where $$A=\frac{R[X_1,\ldots,X_n]^{\dagger}}{(g_1,\ldots,g_r)}$$ with polynomials $g_j\in R[X_1,\ldots,X_n]-\pi R[X_1,\ldots,X_n]$ such that $A$ is $R$-flat. By \cite{vdP} we can lift the $q$-th power Frobenius endomorphism of $\overline{A}$ to an $R$-algebra endomorphism $\sigma$ of $A$. Then $A$, viewed as a $\sigma(A)$-module, is locally free of rank $q^d$. Shrinking $X$ if necessary we may assume that $A$ is a finite free $\sigma(A)$-module of rank $q^d$. As before, $B_K$ denotes a reduced $K$-affinoid algebra, and $B=(B_K)^0$.\\

\addtocounter{satz}{1}{\bf \arabic{section}.\arabic{satz}}\newcounter{matop1}\newcounter{matop2}\setcounter{matop1}{\value{section}}\setcounter{matop2}{\value{satz}} Let $I$ be a countable set. An $I\times I$-matrix ${\cal M}=(a_{i_1,i_2})_{i_1,i_2\in I}$ with entries in an $R$-module $E$ with $E\ne\pi E$ is called {\it nuclear} if for each $M>0$ there are only finitely many $i_2$ such that $\inf_{i_1}\ord_{\pi}({a}_{i_1,i_2})<M$ (thus ${\cal M}$ is nuclear precisely if its transpose is the matrix of a completely continuous operator, or in the terminology of other authors (e.g. \cite{fc3}): a compact operator). An $I\times I$-matrix ${\cal M}=(a_{i_1,i_2})_{i_1,i_2\in I}$ with entries in $A\widehat{\otimes}_RB$ is called {\it nuclear overconvergent} if there exist positive integers $m,c$ and a nuclear matrix $I\times I$-matrix $\widetilde{\cal M}$ with entries in $[m,c]_{B}$ which maps (coefficient-wise) to ${\cal M}$ under the canonical map$$[m,c]_{B}\hookrightarrow R[X]^{\dagger}\widehat{\otimes}_RB\to A\widehat{\otimes}_RB.$$ Clearly, if ${\cal M}$ is nuclear overconvergent then it is nuclear.\\Example: Let $B_K=K$. Nuclear overconvergence implies that the matrix entries are in the subring $A$ of its completion $A\widehat{\otimes}_RR=\widehat{A}$. Conversely, {\it if $I$ is finite}, an $I\times I$-matrix with entries in $A$ is automatically nuclear overconvergent. Similarly, {\it if $I$ is finite}, any $I\times I$-matrix with entries in $\widehat{A}$ is automatically nuclear.\\

\addtocounter{satz}{1}{\bf \arabic{section}.\arabic{satz}} For nuclear matrices ${\cal N}=(c_{h_1,h_2})_{h_1,h_2\in H}$ and ${\cal N}'=(d_{g_1,g_2})_{g_1,g_2\in G}$ with entries in $A\widehat{\otimes}_RB$ define the $(G\times H)\times (G\times H)$-matrix$${\cal N}\otimes{\cal N}':=(e_{(h_1,g_1),(h_2,g_2)})_{(h_1,g_1),(h_2,g_2)\in(G\times H)},$$$$e_{(h_1,g_1),(h_2,g_2)}:=c_{h_1,h_2}d_{g_1,g_2}.$$ Now choose an ordering of the index set $H$. For $k\in\mathbb{N}_0$ let $\bigwedge^k(H)$ be the set of $k$-tuples $(h_1,\ldots,h_k)\in H^k$ with $h_1<\ldots<h_k$. Define the $\bigwedge^k(H)\times\bigwedge^k(H)$-matrix $$\bigwedge^k({\cal N}):={\cal N}^{\wedge k}:=(f_{\vec{h}_1,\vec{h}_2})_{\vec{h}_1,\vec{h}_2\in\bigwedge^k(H)},$$$$f_{\vec{h}_1,\vec{h}_2}=f_{(h_{11},\ldots,h_{1k}),(h_{21},\ldots,h_{2k})}:=\prod_{i=1}^kc_{h_{1,i}h_{2,i}}.$$It is straightforward to check that ${\cal N}\otimes{\cal N}'$ and $\bigwedge^k({\cal N})$ are again nuclear, and even nuclear overconvergent if ${\cal N}$ and ${\cal N'}$ are nuclear overconvergent.\\

\addtocounter{satz}{1}{\bf \arabic{section}.\arabic{satz}}\newcounter{nucmode1}\newcounter{nucmode2}\setcounter{nucmode1}{\value{section}}\setcounter{nucmode2}{\value{satz}} We will use the term "nuclear" also for another concept. Namely, suppose $\psi$ is an operator on a vector space $V$ over $K$. For $g=g(X)\in K[X]$ let $$F(g):=\cap_ng(\psi)^nV\quad\mbox{ and }\quad N(g):=\cup_n\ker g(\psi)^n.$$Let us call a subset $S$ of $K[X]$ bounded away from $0$ if there is an $r\in\mathbb{Q}$ such that $g(a)\ne 0$ for all $\{a\in\mathbb{C}_p;\ord_p(a)\ge r\}$. We say $\psi$ is nuclear if for any subset $S$ of $K[X]$ bounded away from $0$ the following two conditions hold:\\(i) $F(g)\oplus N(g)=V$ for all $g\in S$\\(ii) $N(S):=\sum_{g\in S}N(g)$ is finite dimensional.\\(In particular, if $g\notin(X)$, we can take $S=\{g\}$ and as a consequence of (ii) get $N(g)=\ker g(\psi)^n$ for some $n$.) Suppose $\psi$ is nuclear. Then we can define $P_S(X)=\det(1-X\psi|_{N(S)})$ for subsets $S$ of $K[X]$ bounded away from $0$. These $S$ from a directed set under inclusion, and in \cite{fc3} it is shown that $$P(X):=\lim_{S}P_S(X)$$(coefficient-wise convergence) exists in $K[[X]]$: the characteristic series of $\psi$.\\

\addtocounter{satz}{1}{\bf \arabic{section}.\arabic{satz}} Let $(N_c)_{c\in\mathbb{N}}$ be an inductive system of $B_K$-Banach modules with injective (but not necessarily isometric) transition maps $\rho_{c,c'}:N_c\to N_{c'}$ for $c'\ge c$. Suppose this system has a countable common orthogonal $B_K$-basis, i.e. there is a subset $\{q_m; m\in \mathbb{N}\}$ of $N_1$ such that for all $c$ and $m\in\mathbb{N}$ there are $\lambda_{m,c}\in K^{\times}$ such that $\{\lambda_{m,c}\rho_{1,c}(q_m); m\in\mathbb{N}\}$ is an orthonormal $B_K$-basis of $N_c$. Let $$N:=\lim_{\stackrel{\to}{c}}N_c"="\bigcup_{{c}}N_c$$ and let $N'\subset N$ be a $B_K$-submodule such that $N_c'=N'\cap N_c$ is closed in $N_c$ for all $c$. Endow $N_c'$ with the norm induced from $N_c$ and suppose that also the inductive system $(N_c')_{c\in\mathbb{N}}$ has a countable common orthogonal $B_K$-basis. Let $u$ be a $B_K$-linear endomorphism of $N$ with $u(N')\subset N'$ and restricting to a completely continuous endomorphism $u:N_c\to N_c$ for each $c$. In that situation we have:
 
\begin{pro}\label{quotnuc} $u$ induces a completely continuous $B_K$-endomorphism $u$ of $N_c''=N_c/N_c'$ for each $c$, and $\det(1-uT;N_c'')$ is independent of $c$. If $B_K=K$, the induced endomorphism $u$ of $N''=N/N'$ is nuclear in the sense of \arabic{nucmode1}.\arabic{nucmode2}, and its characteristic series coincides with $\det(1-uT;N_c'')$ for each $c$.
\end{pro}

{\sc Proof:} From \cite{cobmf} A2.6.2 we get that $u$ on $N_c'$ and $u$ on $N_c''$ are completely continuous (note that $N_c''$ is orthonormizable, as follows from \cite{cobmf} A1.2), and that$$\det(1-uT;N_c)=\det(1-uT;N'_c)\det(1-uT;N''_c)$$for each $c$. The assumption on the existence of common orthogonal bases implies (use \cite{eigenc} 4.3.2) $$\det(1-uT;N_c)=\det(1-uT;N_{c'}),\quad\quad \det(1-uT;N'_c)=\det(1-uT;N'_{c'})$$for all $c,c'$. Hence$$\det(1-uT;N''_c)=\det(1-uT;N''_{c'})$$for all $c,c'$. Also note that for $c'\ge c$ the maps $N_c''\to N_{c'}''$ are injective. The additional assumptions in case $B_K=K$ now follow from \cite{fc3} Theorem 1.3 and Lemma 1.6.\\

\addtocounter{satz}{1}{\bf \arabic{section}.\arabic{satz}}\newcounter{dwode1}\newcounter{dwode2}\setcounter{dwode1}{\value{section}}\setcounter{dwode2}{\value{satz}} Shrinking $X$ if necessary we may assume that the module of ($p$-adically separated) differentials $\Omega^1_{A/R}$ is free over $A$. Fix a basis $\omega_1,\ldots,\omega_d$. With respect to this basis, let ${\cal D}$ be the $d\times d$-matrix of the $\sigma$-linear endomorphism of $\Omega^1_{A/R}$ which the $R$-algebra endomorphism $\sigma$ of $A$ induces. Then ${\cal D}^{\wedge k}=\bigwedge^k({\cal D})$ is the matrix of the $\sigma$-linear endomorphism of $\Omega^k_{A/R}=\bigwedge^k(\Omega^1_{A/R})$ which $\sigma$ induces.\\Let $\theta=\sigma^{-1}\circ \tr$ be the endomorphism of $\Omega^d_{A/R}$ constructed in \cite{fc1} Theorem 8.5. It is a Dwork operator: we have $\theta(\sigma(a)y)=a\theta(y)$ for all $a\in A$, $y\in \Omega^d_{A/R}$. Denote also by $\theta$ the Dwork operator on $A$ which we get by transport of structure from $\theta$ on $\Omega^d_{A/R}$ via the isomorphism $A\cong\Omega^d_{A/R}$ which sends $1\in A$ to our distinguished basis element $\omega_{1}\wedge\ldots\wedge\omega_{d}$ of $\Omega^d_{A/R}$.\\For $c\in\mathbb{N}$ define the subring $A^c$ of $A_K=A\otimes_RK$ as the image of $$T_n^c\hookrightarrow R[X]^{\dagger}\otimes_RK\to A_K.$$This is again a $K$-affinoid algebra, and we have \label{randink} $$\theta(A^c)\subset A^c$$ for $c>>0$. To see this, choose an $R$-algebra endomorphism $\widetilde{\sigma}$ of $R[X]^{\dagger}$ which lifts both $\sigma$ on $A$ and the $q$-th power Frobenius endomorphism on $k[X]$. With respect to this $\widetilde{\sigma}$ choose a Dwork operator $\widetilde{\theta}$ on $R[X]^{\dagger}$ lifting $\theta$ on $A$ (as in the beginning of the proof of \cite{fc3} Theorem 2.3). Then apply \cite{fc3} Lemma 2.4 which says $\widetilde{\theta}(T_n^c)\subset T_n^c$.\\

\addtocounter{satz}{1}{\bf \arabic{section}.\arabic{satz}}\newcounter{dwrat1}\newcounter{dwrat2}\setcounter{dwrat1}{\value{section}}\setcounter{dwrat2}{\value{satz}} Let ${\cal M}=(a_{i_1,i_2})_{i_1,i_2\in I}$ be a nuclear overconvergent $I\times I$-matrix with entries in $A\widehat{\otimes}_RB$. For $c\in\mathbb{N}$ let $\check{M}_I^c$ be the $A^c\widehat{\otimes}_KB_K$-Banach module for which the set of symbols $\{\check{e}_i\}_{i\in I}$ is an orthonormal basis. For $c\ge c'$ we have the continuous inclusion of $B_K$-algebras $A^{c'}\widehat{\otimes}_KB_K\subset A^c\widehat{\otimes}_KB_K$, hence a continuous inclusion of $B_K$-modules $\check{M}_I^{c'}\subset\check{M}_I^c$. Since ${\cal M}$ is nuclear overconvergent we have $a_{i_1,i_2}\in A^c\widehat{\otimes}_KB_K$ for all $c>>0$, all $i_1, i_2$. We may thus define for all $c>>0$ the $B_K$-linear endomorphism $\psi=\psi[{\cal M}]$ of $\check{M}_I^c$ by$$\psi(\sum_{i_1\in I}b_{i_1}\check{e}_{i_1})=\sum_{i_1\in I}\sum_{i_2\in I}(\theta\otimes1)(b_{i_1}a_{i_1,i_2})\check{e}_{i_2}$$($b_{i_1}\in A^c\widehat{\otimes}_KB_K$). Clearly these endomorphisms extend each other for increasing $c$, hence we get an endomorphism $\psi=\psi[{\cal M}]$ on $$\check{M}_I:=\bigcup_{c>>0}\check{M}^c_I.$$

\addtocounter{satz}{1}{\bf \arabic{section}.\arabic{satz}}\newcounter{dwint1}\newcounter{dwint2}\setcounter{dwint1}{\value{section}}\setcounter{dwint2}{\value{satz}} Suppose $B_K=K$ and $I$ is finite, and ${\cal M}$ is the matrix of the $\sigma$-linear endomorphism $\phi$ acting on the basis $\{e_i\}_{i\in I}$ of the free $A$-module $M$. Then we define $\psi[{\cal M}]$ as the Dwork operator $$\psi[{\cal M}]:\ho_{A}(M,\Omega^d_{A/R})\to\ho_{A}(M,\Omega^d_{A/R}),\quad f\mapsto \theta\circ f\circ\phi.$$This definition is compatible with that in \arabic{dwrat1}.\arabic{dwrat2}: Consider the canonical embedding$$\ho_{A}(M,\Omega^d_{A/R})\to \ho_{A}(M,\Omega^d_{A/R})\otimes_RK\stackrel{w}{\cong}\check{M}_I$$where the  inverse of the $A_K$-linear isomorphism $w$ sends $\check{e}_i\in\check{M}_I$ to the homomorphism which maps $e_i\in M$ to $\omega_{1}\wedge\ldots\wedge\omega_{d}$ and which maps $e_{i'}$ for $i'\ne i$ to $0$. This embedding commutes with the operators $\psi[{\cal M}]$.\\

\begin{satz}\label{nucl} For each $c>>0$, the endomorphism $\psi=\psi[{\cal M}]$ on $\check{M}^c_I$ is a completely continuous $B_K$-Banach module endomorphism. Its Fredholm determinant $\det(1-\psi T;\check{M}^c_I)$ is independent of $c$. Denote it by $\det(1-\psi T;\check{M}_I)$. If $B_K=K$, the endomorphism $\psi=\psi[{\cal M}]$ on $\check{M}_I$ is nuclear in the sense of \cite{fc3}, and its characteristic series as defined in \cite{fc3} coincides with $\det(1-\psi T;\check{M}_I)$.
\end{satz}

{\sc Proof:} Choose a lifting of ${\cal M}=(a_{i_1,i_2})_{i_1,i_2\in I}$ to a nuclear matrix $(\widetilde{a}_{i_1,i_2})_{i_1,i_2\in I}$ with entries in $[m,c]_{B}$. Also choose a lifting of $\theta$ on $A$ to a Dwork operator $\widetilde{\theta}$ on $R[X]^{\dagger}$ (with respect to a lifting of $\sigma$, as in \arabic{dwode1}.\arabic{dwode2}). Let $N_I^c$ be the $T_n^c\widehat{\otimes}_KB_K$-Banach module for which the set of symbols $\{(\check{e}_i)\widetilde{}\}_{i\in I}$ is an orthonormal basis, and define the $B_K$-linear endomorphism $\widetilde{\psi}$ of $N_I^c$ by$$\widetilde{\psi}(\sum_{i_1\in I}\widetilde{b}_{i_1}(\check{e}_{i_1})\widetilde{})=\sum_{i_1\in I}\sum_{i_2\in I}(\widetilde{\theta}\otimes1)(\widetilde{b}_{i_1}\widetilde{a}_{i_1,i_2})(\check{e}_{i_2})\widetilde{}$$($\widetilde{b}_{i_1}\in T_n^c\widehat{\otimes}_KB$). An orthonormal basis of $N_I^c$ as a $B_K$-Banach module is given by \begin{gather}\{\pi^{[\frac{|\alpha|}{c}]}X^{\alpha}(\check{e}_i)\widetilde{}\}_{\alpha\in\mathbb{N}_0^n, i\in I}.\tag{$1$}\end{gather}By \ref{matlim} the matrix for $\widetilde{\psi}$ in this basis is completely continuous; that is, $\widetilde{\psi}$ is completely continuous. If $I_c\subset T_n^c$ and $I_{\infty}\subset T_n$ denote the respective ideals generated by the elements $g_1,\ldots,g_r$ from \arabic{setgeo1}.\arabic{setgeo2}, then $I_{\infty}\cap T_n^c$ is the kernel of $T_n^c\to A^c$, so by \ref{exaide} the sequences \begin{gather}0\to I_c\to T_n^c\to A^c\to0\tag{$2$}\end{gather}are exact for $c>>0$. Let $H$ be the $B_K$-Banach module with orthonormal basis the set of symbols $\{h_i\}_{i\in I}$. From $(2)$ we derive an exact sequence\begin{gather}0\to I_c\widehat{\otimes}_KH\to T_n^c\widehat{\otimes}_KH\to A^c\widehat{\otimes}_KH\to0\tag{$3$}\end{gather}(To see exactness of $(3)$ on the right note that one of the equivalent norms on $A^c$ is the residue norm for the surjective map of $K$-affinoid algebras $T_n^c\to A^c$ (this surjection even has a continuous $K$-linear section as the proof of \cite{cobmf} A2.6.2 shows)). We use the following isomorphisms of $T_n^c\widehat{\otimes}_KB_K$-Banach modules (in $(i)$) resp. of $A^c\widehat{\otimes}_KB_K$-Banach modules (in $(ii)$):\begin{gather}T_n^c\widehat{\otimes}_KH=(T_n^c\widehat{\otimes}_KB_K)\widehat{\otimes}_{B_K}H\cong N_I^c,\quad 1\otimes h_i\mapsto(\check{e}_i)\widetilde{}\tag{$i$}\end{gather} \begin{gather}A^c\widehat{\otimes}_KH=(A^c\widehat{\otimes}_KB_K)\widehat{\otimes}_{B_K}H\cong\check{M}_I^c,\quad1\otimes h_i\mapsto\check{e}_i\tag{$ii$}\end{gather} By \ref{onbas} we find a subset $E$ of $\mathbb{N}_0^n\times\{1,\ldots,r\}$ such that $\{\pi^{[\frac{|\alpha|+d_j}{c}]}X^{\alpha}g_j\}_{(\alpha,j)\in E}$ is an orthonormal basis of $I_{c}$ over $K$ for all $c>>0$. For the $B_K$-Banach modules $I_c\widehat{\otimes}_KH=(I_c\widehat{\otimes}_KB_K)\widehat{\otimes}_{B_K}H$ we therefore have the orthonormal basis\begin{gather}\{\pi^{[\frac{|\alpha|+d_j}{c}]}X^{\alpha}g_j\otimes h_i\}_{(\alpha,j)\in E, i\in I}.\tag{$4$}\end{gather} It is clear that the systems of orthonormal bases $(1)$ resp. $(4)$ make up systems of common orthonormal bases when $c$ increases. (This is why we took pains to prove \ref{onbas}; the present argument could be simplified if we could prove the existence of a common orthogonal basis for the system $(\check{M}^c_I)_{c>>0}$.) Now let $N_I=\cup_{c}N_I^c$. From the exactness of the sequences $(3)$ and from the injectivity of the maps $\check{M}^c_I\to \check{M}^{c'}_I$ for $c\le c'$ we get $I_c\widehat{\otimes}_KH=T_n^c\widehat{\otimes}_KH\cap\ke(N_I\to\check{M}_I)$. Thus the theorem follows from \ref{quotnuc}.\\

\begin{kor}\label{fredganz} $$\prod_{r=0}^d\det(1-\psi[{\cal M}\otimes{\cal D}^{\wedge(d-r)}] T;\check{M}_I)^{(-1)^{r-1}}$$is the quotient of entire power series in the variable $T$ with coefficients in $B_K$; in other words, it is a meromorphic function on $\mathbb{A}^1_K\times_{\spm(K)}\spm(B_K)$.\\
\end{kor}

\addtocounter{satz}{1}{\bf \arabic{section}.\arabic{satz}}\newcounter{lfnd1}\newcounter{lfnd2}\setcounter{lfnd1}{\value{section}}\setcounter{lfnd2}{\value{satz}} Let $B_K=K$. We want to define the $L$-function of a nuclear matrix ${\cal M}=(a_{i_1,i_2})_{i_1,i_2\in I}$ (with entries in $\widehat{A}$). For $f\in\mathbb{N}$ define the $f$-fold $\sigma$-power ${\cal M}^{(\sigma)^f}$ of ${\cal M}$ to be the matrix product $${\cal M}^{(\sigma)^f}:=((a_{i_1,i_2})_{i_1,i_2\in I})(\sigma(a_{i_1,i_2})_{i_1,i_2\in I})\ldots((\sigma^{f-1}(a_{i_1,i_2}))_{i_1,i_2\in I}).$$Let $\overline{x}\in X$ be a geometric point of degree $f$ over $k$, that is, a surjective $k$-algebra homomorphism $\overline{A}\to\mathbb{F}_{q^f}$. Let $R_f$ be the unramified extension of $R$ with residue field $\mathbb{F}_{q^f}$, and let $x:\widehat{A}\to R_f$ be the Teichm\"uller lifting of $\overline{x}$ with respect to $\sigma$ (the unique $\sigma^f$-invariant surjective $R$-algebra homomorphism lifting $\overline{x}$). By (quite severe) abuse of notation we write$${\cal M}_{\overline{x}}:=x({\cal M}^{(\sigma)^f}),$$the $I\times I$-matrix with $R_f$-entries obtained by applying $x$ to the entries of ${\cal M}^{(\sigma)^f}$ --- the "fibre of ${\cal M}$ in ${\overline{x}}$". The nuclearity condition implies that ${\cal M}_{\overline{x}}$ is nuclear; equivalently, its transpose is a completely continuous matrix over $R_f$ in the sense of \cite{ser}. It turns out that the Fredholm determinant $\det(1-{\cal M}_{\overline{x}}T^{\deg(\overline{x})})$ has coefficients in $R$, not just in $R_f$. We define the $R[[T]]$-element$$L({\cal M},T):=\prod_{\overline{x}\in X}\frac{1}{\det(1-{\cal M}_{\overline{x}}T^{\deg(\overline{x})})}.$$It is trivially holomorphic on the open unit disk. Let $T$ be the set of $k$-valued points $\overline{x}:\overline{A}\to k$ of $X$. For a completely continuous endomorphism ${\psi}$ of an orthonormizable $K$-Banach module we denote by $\tr_K({\psi})\in K$ its trace.\\

\begin{satz}\label{tracf} Let ${\cal M}$ be a nuclear overconvergent matrix over $\hat{A}$.\\(1) For each $\overline{x}\in T$ the element $$S_{\overline{x}}:=\sum_{0\le j\le d}(-1)^j\tr(({\cal D}^{\wedge d-j})_{\overline{x}})$$is invertible in $R$. For $0\le i\le d$, we have$$\tr_K(\psi[{\cal M}\otimes{\cal D}^{\wedge d-i}])=\sum_{\overline{x}\in T}\frac{\tr(({\cal D}^{\wedge d-i})_{\overline{x}})\tr({\cal M}_{\overline{x}})}{S_{\overline{x}}}.$$(2) $${L}({\cal M},T)=\prod_{r=0}^d\det(1-\psi[{\cal M}\otimes{\cal D}^{\wedge(d-r)}] T;\check{M}_I)^{(-1)^{r-1}}.$$In particular, by \ref{fredganz}, ${L}({\cal M},T)$ is meromorphic on $\mathbb{A}_K^1$.\\
\end{satz}

{\sc Proof:} Let $J\subset A$ be the ideal generated by all elements of the form $a-\sigma(a)$ with $a\in A$. Then $\spec(A/J)$ is a direct sum of copies of $\spec(R)$, indexed by $T$: It is the direct sum of all Teichm\"uller lifts of elements in $T$ (or rather, their restrictions from $\widehat{A}$ to $A$; cf. \cite{fc3} Lemma 3.3). Let $C(A,\sigma)$ be the category of finite (not necessarily projective) $A$-modules $(M,\phi)$ with a $\sigma$-linear endomorphism $\phi$, let $m(A,\sigma)$ be the free abelian group generated by the isomorphism classes of objects of $C(A,\sigma)$, and let $n(A,\sigma)$ be the subgroup of $m(A,\sigma)$ generated by the following two types of elements. The first type is of the form $(M,\phi)-(M_1,\phi_1)-(M_2,\phi_2)$ where$$0\to(M_1,\phi_1)\to(M,\phi)\to(M_2,\phi_2)\to0$$is an exact sequence in $C(A,\sigma)$. The second type is of the form $(M,\phi_1+\phi_2)-(M,\phi_1)-(M,\phi_2)$ for $\sigma$-linear operators $\phi_1$, $\phi_2$ on the same $M$. Set $K(A,\sigma)=m(A,\sigma)/n(A,\sigma)$. By the analogous procedure define the group $K^*(A,\sigma)$ associated  with the category of finite $A$-modules with a Dwork operator relative to $\sigma$. (Here we follow the notation in \cite{wahrk}. The notation in \cite{fc3} is the opposite one !). By \cite{fc3}, both $K(A,\sigma)$ and $K^*(A,\sigma)$ are free $A/J$-modules of rank one. For a finite square matrix ${\cal N}$ over $A$ we denote by $\tr_{A/J}({\cal N})\in A/J$ the trace of the matrix obtained by reducing modulo $J$ the entries of ${\cal N}$. Moreover, for such ${\cal N}$ we view $\psi[{\cal N}]$ always as a Dwork operator on a (finite) $A$-module as in \arabic{dwint1}.\arabic{dwint2}, i.e. we do not invert $\pi$. From \cite{wahrk} sect.3 it follows that $\psi[{\cal D}^{\wedge d-i}]$ can be identified with the standard Dwork operator $\psi_i$ on $\Omega^i_{A/R}$ from \cite{fc3}. By \cite{fc3} sect.5 Cor.1 we have \begin{gather}[\psi[{\cal D}^{\wedge 0}]]\sum_{0\le j\le d}(-1)^j\tr_{A/J}({\cal D}^{\wedge d-j})=[(A,\id)]\tag{$1$}\end{gather}in $K^*(A,\sigma)$, and $\sum_{0\le j\le d}(-1)^j\tr_{A/J}({\cal D}^{\wedge d-j})$ is invertible in $A/J$. By \cite{fc3} Theorem 5.2 we also have \begin{gather}[\psi[{\cal D}^{\wedge d-i}]]=\tr_{A/J}({\cal D}^{\wedge d-i})[\psi[{\cal D}^{\wedge 0}]]\tag{$2$}\end{gather}in $K^*(A,\sigma)$. To prove the theorem suppose first that ${\cal M}$ is a finite square matrix. It then gives rise to an element $[{\cal M}]$ of $K(A,\sigma)$. By \cite{wahrk} 10.8 we have $$[{\cal M}]=\tr_{A/J}({\cal M})[(A,\id)]$$in $K(A,\sigma)$. Application of the homomorphism of $A/J$-modules $$\lambda_{i}:K(A,\sigma)\to K^*(A,\sigma)$$of \cite{wahrk} p.42 gives \begin{gather}[\psi[{\cal M}\otimes{\cal D}^{\wedge d-i}]]=\tr_{A/J}({\cal M})[\psi[{\cal D}^{\wedge d-i}]]\tag{$3$}\end{gather} in $K^*(A,\sigma)$. From $(1), (2), (3)$ we get $$[\psi[{\cal M}\otimes{\cal D}^{\wedge d-i}]]=\frac{\tr_{A/J}({\cal M})\tr_{A/J}({\cal D}^{\wedge d-i})}{\sum_{0\le j\le d}(-1)^j\tr_{A/J}({\cal D}^{\wedge d-j})}[(A,\id)]$$in $K^*(A,\sigma)$. Taking the $R$-trace proves (1) in case ${\cal M}$ is a finite square matrix. Then taking the alternating sum over $0\le i\le d$ gives the additive formulation of (2) in case ${\cal M}$ is a finite square matrix (see also \cite{wahrk} Theorem 3.1).\\The case where the index set $I$ for ${\cal M}$ is infinite follows by a limiting argument from the case where $I$ is finite. We explain this for (2), leaving the easier (1) to the reader. Let ${\cal P}(I)$ be the set of {\it finite} subsets of $I$. For $I'\in{\cal P}(I)$, the $I'\times I'$-sub-matrix ${\cal M}^{I'}=(a_{i_1,i_2})_{i_1,i_2\in I'}$ of ${\cal M}$ is again nuclear overconvergent. Hence, in view of the finite square matrix case it is enough to show \begin{gather}{L}({\cal M},T)=\lim_{I'\in{\cal P}(I)}{L}({\cal M}^{I'},T)\tag{$1$}\end{gather}and for any $0\le r\le d$ also \begin{gather}\det(1-\psi[{\cal M}\otimes{\cal D}^{\wedge r}] T;\check{M}_I)=\lim_{I'\in{\cal P}(I)}\det(1-\psi[{\cal M}^{I'}\otimes{\cal D}^{\wedge r}] T;\check{M}_I)\tag{$2$}\end{gather}(coefficient-wise convergence). For $I'\in{\cal P}(I)$ define the $I\times I$-matrix ${\cal M}[I']=(a^{I'}_{i_1,i_2})_{i_1,i_2\in I}$ by $a^{I'}_{i_1,i_2}=a_{i_1,i_2}$ if $i_2\in I'$ and $a^{I'}_{i_1,i_2}=0$ otherwise. For a geometric point $\overline{x}\in X$ we may view the fibre matrices ${\cal M}_{\overline{x}}$ resp. ${\cal M}[I']_{\overline{x}}$ for $I'\in{\cal P}(I)$ as the transposed matrices of completely continuous operators $\lambda_{\overline{x}}$ resp. $\lambda[I']_{\overline{x}}$ acting all on one single $K$-Banach space $E_{\overline{x}}$ with orthonormal basis indexed by $I$. And we may view the fibre matrix ${\cal M}^{I'}_{\overline{x}}$ as the transposed matrix of the restriction of $\lambda[I']_{\overline{x}}$ to a $\lambda[I']_{\overline{x}}$-stable subspace of $E_{\overline{x}}$, spanned by a finite subset of our given orthonormal basis and containing $\lambda[I']_{\overline{x}}(E_{\overline{x}})$. For the norm topology on the space $L(E_{\overline{x}},E_{\overline{x}})$ of continuous $K$-linear endomorphisms of $E_{\overline{x}}$ we find, using the nuclearity of ${\cal M}$, that $\lim_{I'}\lambda[I']_{\overline{x}}=\lambda_{\overline{x}}$. Hence it follows from \cite{ser} prop.7,c) that$$\det(1-{\cal M}_{\overline{x}}T^{\deg(\overline{x})})=\lim_{I'\in{\cal P}(I)}\det(1-{\cal M}[I']_{\overline{x}}T^{\deg(\overline{x})}).$$But by \cite{ser} prop.7,d) we have $$\det(1-{\cal M}[I']_{\overline{x}}T^{\deg(\overline{x})})=\det(1-{\cal M}^{I'}_{\overline{x}}T^{\deg(\overline{x})}).$$Together we get $(1)$. The proof of $(2)$ is similar: By the proof of \ref{matlim} we have indeed$$\lim_{I'\in{\cal P}(I)}\psi[{\cal M}[I']\otimes{\cal D}^{\wedge r}]=\psi[{\cal M}\otimes{\cal D}^{\wedge r}]$$in the space of continuous $K$-linear endomorphisms of $\check{M}^c_I$, so \cite{ser} prop.7,c) gives$$\det(1-\psi[{\cal M}\otimes{\cal D}^{\wedge r}] T;\check{M}^c_I)=\lim_{I'\in{\cal P}(I)}\det(1-\psi[{\cal M}[I']\otimes{\cal D}^{\wedge r}] T;\check{M}^c_I).$$ Now the $\psi[{\cal M}[I']\otimes{\cal D}^{\wedge r}]$ do not have finite dimensional image in general, but clearly an obvious generalization of \cite{ser} prop.7,d) shows$$\det(1-\psi[{\cal M}[I']\otimes{\cal D}^{\wedge r}] T;\check{M}^c_I)=\det(1-\psi[{\cal M}^{I'}\otimes{\cal D}^{\wedge r}] T;\check{M}^c_I)$$for $I'\in{\cal P}(I)$. We are done.\\  

\section{The Grothendieck group}

In this section we introduce the Grothendieck group $\Delta(A\widehat{\otimes}_RB)$ of nuclear $\sigma$-modules. It is useful since on the one hand, formation of the $L$-function of a given nuclear $\sigma$-module factors over this group, and on the other hand, many natural nuclear $\sigma$-modules which are not nuclear overconvergent can be represented in this group through nuclear overconvergent ones.\\ 

\addtocounter{satz}{1}{\bf \arabic{section}.\arabic{satz}} We will write $\sigma$ also for the endomorphism $\sigma\otimes1$ of $A\widehat{\otimes}_RB=\widehat{A}\widehat{\otimes}_RB$. For $\ell=1,2$ let ${\cal M}_{\ell}$ be $I_{\ell}\times I_{\ell}$-matrices with entries in $A\widehat{\otimes}_RB$, for countable index sets $I_{\ell}$. We say ${\cal M}_1$ is {\it $\sigma$-similar} to ${\cal M}_2$ over $A\widehat{\otimes}_RB$ if there exist a $I_1\times I_2$-matrix ${\cal S}$ and a $I_2\times I_1$-matrix ${\cal S}'$, both with entries in $A\widehat{\otimes}_RB$, such that ${\cal S}{\cal S}'$ (resp. ${\cal S}'{\cal S}$) is the identity $I_1\times I_1$ (resp. $I_2\times I_2$) -matrix, and such that ${\cal S}'{\cal M}_1{\cal S}^{\sigma}={\cal M}_2$ (in particular it is required that all these matrix products converge coefficient-wise in $A\widehat{\otimes}_RB$). Clearly, $\sigma$-similarity is an equivalence relation.\\

\addtocounter{satz}{1}{\bf \arabic{section}.\arabic{satz}} Let $m(A\widehat{\otimes}_RB)$ be the free abelian group generated by the $\sigma$-similarity classes of nuclear matrices (over arbitrary countable index sets) with entries in $A\widehat{\otimes}_RB$. Let $\Delta(A\widehat{\otimes}_RB)$ be the quotient of $m(A\widehat{\otimes}_RB)$ by the subgroup generated by all the elements $[{\cal M}]-[{\cal M}']-[{\cal M}'']$ for matrices ${\cal M}=(a_{i_1,i_2})_{i_1,i_2\in I}$, ${\cal M}'=(a_{i_1,i_2})_{i_1,i_2\in I'}$ and ${\cal M}''=(a_{i_1,i_2})_{i_1,i_2\in I''}$ where $I=I'\sqcup I''$ is a partition of $I$ such that $a_{i_1,i_2}=0$ for all pairs $(i_1,i_2)\in I'\times I''$ (in other words, ${\cal M}$ is in block triangular form and ${\cal M}'$, ${\cal M}''$ are the matrices on the block diagonal). 

Elements $z\in \Delta(A\widehat{\otimes}_RB)$ can be written as $z=[{\cal M}_{+}]-[{\cal M}_{-}]$ with nuclear matrices ${\cal M}_{+}, {\cal M}_{-}$. If $\{{\cal M}_n\}_{n\in\mathbb{N}}$ is a collection of nuclear matrices such that $\ord_{\pi}({\cal M}_n)\to\infty$ (where $\ord_{\pi}({\cal M})=\min_{i_1,i_2}\{\ord_{\pi}a_{i_1,i_2}\}$ for a matrix ${\cal M}=(a_{i_1,i_2})_{i_1,i_2\in I}$) and if $\{\nu_n\}_{n\in\mathbb{N}}$ are integers, then the infinite sum $\sum_{n\in\mathbb{N}}\nu_n[{\cal M}_n]$ can be viewed as an element of $\Delta(A\widehat{\otimes}_RB)$ as follows: Sorting the $\nu_n$ according to their signs means breaking up this sum into a positive and a negative summand, so we may assume $\nu_n\ge 1$ for all $n$. Replacing ${\cal M}_n$ by the block diagonal matrix $diag({\cal M}_n,{\cal M}_n,\ldots,{\cal M}_n)$ with $\nu_n$ copies of ${\cal M}_n$ we may assume $\nu_n=1$ for all $n$. Since all ${\cal M}_n$ are nuclear and $\ord_{\pi}({\cal M}_n)\to\infty$ the block diagonal matrix ${\cal M}=diag({\cal M}_1,{\cal M}_2,{\cal M}_3,\ldots)$ is nuclear. It represents the desired element of $\Delta(A\widehat{\otimes}_RB)$. Matrix tensor product (see \arabic{matop1}.\arabic{matop2}) defines a multiplication in $\Delta(A\widehat{\otimes}_RB)$: One checks that$$([{\cal M}_{1,+}]-[{\cal M}_{1,-}])\otimes([{\cal M}_{2,+}]-[{\cal M}_{2,-}])$$$$=[{\cal M}_{1,+}\otimes{\cal M}_{2,+}]-[{\cal M}_{1,+}\otimes{\cal M}_{2,-}]-[{\cal M}_{1,-}\otimes{\cal M}_{2,+}]+[{\cal M}_{1,-}\otimes{\cal M}_{2,-}]$$is independent of the chosen representations.\\  

\addtocounter{satz}{1}{\bf \arabic{section}.\arabic{satz}} A more suggestive way to think of $\Delta(A\widehat{\otimes}_RB)$ is the following. We say that a subset $\{e_i\}_{i\in I}$ of an $A\widehat{\otimes}_RB$-module $M$ is a {\it formal basis} if there is an isomorphism of $A\widehat{\otimes}_RB$-modules $$\{(d_i)_{i\in I};\, d_i\in A\widehat{\otimes}_RB\}\cong M$$mapping for any $j\in I$ the sequence $(d_i)_i$ with $d_j=1$ and $d_i=0$ for $i\ne j$ to $e_j$. A {\it nuclear} $\sigma$-{\it module over} $A\widehat{\otimes}_RB$ is an $A\widehat{\otimes}_RB$-module $M$ together with a $\sigma$-linear endomorphism $\phi$ such that there exists a formal basis $\{e_i\}_{i\in I}$ of $M$ such that the action of $\phi$ on $\{e_i\}_{i\in I}$ is described by a nuclear matrix ${\cal M}$ with entries in $A\widehat{\otimes}_RB$, i.e. $\phi e_i={\cal M}e_i$ if we think of $e_i$ as the $i$-th column of the identity $I\times I$ matrix. We usually think of a nuclear $\sigma$-module over $A\widehat{\otimes}_RB$ as a {\it family of nuclear $\sigma$-modules over $A$, parametrized by the rigid space $\spm(B_K)$}.\\In the above situation, if ${\cal S}$ is a (topologically) invertible $I\times I$-matrix with entries in $A\widehat{\otimes}_RB$, then ${\cal S}^{-1}{\cal M}{\cal S}^{\sigma}$ is the matrix of $\phi$ in the new formal basis consisting of the elements ${\cal S}e_i=e'_i$ of $M$ (if now we think of $e'_i$ as the $i$-th column of the identity $I\times I$-matrix). Hence we can view $\Delta(A\widehat{\otimes}_RB)$ as the Grothendieck group of nuclear $\sigma$-modules over $A\widehat{\otimes}_RB$, i.e. as the quotient of the free abelian group generated by (isomorphism classes of) nuclear $\sigma$-modules over $A\widehat{\otimes}_RB$, divided out by the relations $[(M,\phi)]-[(M',\phi')]-[(M'',\phi'')]$ coming from short exact sequences$$0\to (M',\phi')\to(M,\phi)\to(M'',\phi'')\to0$$which are $A\widehat{\otimes}_RB$-linearly (but not necessarily $\phi$-equivariantly) split.\\

\begin{pro}\label{lfndbkt} Let $B_K=K$. Let $x\in\Delta(\widehat{A})$ be represented by a convergent series $x=\sum_{\ell\in\mathbb{N}}\nu_{\ell}[{\cal M}_{\ell}]$ with nuclear matrices ${\cal M}_{\ell}$ over $\widehat{A}$. Then the $L$-series$$L(x,T):=\prod_{\ell\in\mathbb{N}}L({\cal M}_{\ell},T)^{\nu_{\ell}}$$is independent of the chosen representation of $x$. If all ${\cal M}_{\ell}$ are nuclear overconvergent, then $L(x,T)$ represents a meromorphic function on $\mathbb{A}_K^1$.
\end{pro}

{\sc Proof:} One checks that $\sigma$-similar nuclear matrices over $\widehat{A}$ have the same $L$-function. Indeed, even the Euler factors at closed points of $X$ are the same: they are given by Fredholm determinants of similar (in the ordinary sense) completely continuous matrices. Now let ${\cal M}$, ${\cal M}'$ and ${\cal M}''$ give rise to a typical relation $[{\cal M}]=[{\cal M}']+[{\cal M}'']$ as in our definition of $\Delta(A\widehat{\otimes}_RB)$. Then one checks that $$L({\cal M},T)=L({\cal M}',T)L({\cal M}'',T),$$ again by comparing Euler factors. And finally it also follows from the Euler product definition that $\ord_{\pi}(1-L({\cal M}_{\ell},T))\to \infty$ if $\ord_{\pi}({\cal M}_{\ell})\to \infty$. Altogether we get the well definedness of $L(x,T)$. If the ${\cal M}_l$ are nuclear overconvergent, then the $L({\cal M}_{\ell},T)$ are meromorphic by \ref{tracf} and we get the second assertion.\\

\section{Resolution of unit root parts of rank one}
\label{resol}

In this section we describe a family version of the limiting module construction. Given a rank one unit root $\sigma$-module $(M_{unit},\phi_{unit})$ which is the unit root part of a (unit root ordinary) nuclear $\sigma$-module $(M,\phi)$ and such that $\phi_{unit}$ acts by a $1$-unit $a_{i_0,i_0}\in\widehat{A}$ on a basis element of $M_{unit}$, we choose an affinoid rigid subspace $\spm(B_K)$ of $\mathbb{A}_K^1$ such that for each $\mathbb{C}_p$-valued point $x\in \spm(B_K)\subset\mathbb{C}_p$ the exponentiation $a_{i_0,i_0}^x$ is well defined. Hence we get a rank one $\sigma$-module over $A\widehat{\otimes}_RB$. We express its class in $\Delta(\widehat{A}\widehat{\otimes}_RB)$ through a set (indexed by $r\in\mathbb{Z}$) of nuclear $\sigma$-modules $(B^{r}(M),B^r(\phi))$ over $A\widehat{\otimes}_RB$ which are overconvergent if $(M,\phi)$ is overconvergent, even if $(M_{unit},\phi_{unit})$ is not overconvergent. Later $\spm(B_K)$ will be identified with the set of characters $\kappa:U_{R}^{(1)}\to\mathbb{C}_p^{\times}$ of the type $\kappa(u)=\kappa_x(u)=u^x$ for small $x\in\mathbb{C}_p$, where $U_{R}^{(1)}$ denotes the group of $1$-units in $R$. To obtain the optimal parameter space for the $B^{r}(M)$ (i.e. the maximal region in $\mathbb{C}_p$ of elements $x$ for which $\kappa_x$ occurs in the parameter space) one needs to go to the union of all these $\spm(B_K)$. This $K$-rigid space is not affinoid any more; in the case $K=\mathbb{Q}_p$ it is the parameter space ${\cal B}^*$ from \cite{cobmf}. We will however not pass to this limit here, since for an extension of the associated unit root $L$-function even to the {\it whole} character space we will have another method available in section \ref{mostnovel}.\\ 

\begin{lem}\label{logexp} Let $E$ be a $p$-adically separated and complete ring such that $E\to E\otimes\mathbb{Q}$ is injective and denote again by $\ord_p$ the natural extension of $\ord_p$ from $E$ to $E\otimes\mathbb{Q}$. \\(i) Let $x\in E$. If $\ord_{p}(x)>\frac{1}{p-1}$, then $\ord_{p}(\frac{x^n}{n!})\ge0$ for all $n\ge0$, and $$\exp(x)=\sum_{n\ge0}\frac{x^n}{n!}$$ converges.\\(ii) Let $x\in E$. If $\ord_{p}(x)>0$, then $\ord_p(\frac{x^n}{n})\ge0$ for all $n\ge1$, and $$\log(1+x)=\sum_{n\ge1}(-1)^{n-1}\frac{x^n}{n}$$ converges. Moreover, if $\ord_{p}(x)>\beta\ge \frac{1}{p-1}$, then $\ord_{p}(\log(1+x))>\beta$; if $\ord_{p}(x)\ge \frac{1}{p-1}\frac{1}{p^b}$ for some $b\in\mathbb{N}_0$, then $\ord_{p}(\log(1+x))\ge\frac{1}{p-1}-b$.\end{lem}

{\sc Proof:} Proceed as in \cite{rob}, p.252, p.356.\\

\addtocounter{satz}{1}{\bf \arabic{section}.\arabic{satz}}\newcounter{nordef1}\newcounter{nordef2}\setcounter{nordef1}{\value{section}}\setcounter{nordef2}{\value{satz}} Fix a countable non empty set $I$ and an element $i_0\in I$, let $I_1=I-\{i_0\}$. Let ${\cal M}=(a_{i_1,i_2})_{i_1,i_2\in I}$ be a nuclear $I\times I$-matrix over $\widehat{A}$. It is called $1$-{\it normal} if $1-a_{i_0,i_0}\in\pi \widehat{A}$ and if $a_{i_1,i_2}\in\pi \widehat{A}$ for all $(i_1,i_2)\ne(i_0,i_0)$. It is called {\it standard normal} if $a_{i_1,i_0}=0$ for all $i_1\in I_1$, if $a_{i_0,i_0}$ is invertible in $\widehat{A}$ and if $a_{i_1,i_2}\in\pi \widehat{A}$ for all $(i_1,i_2)\ne(i_0,i_0)$. It is called {\it standard} $1$-{\it normal} if it is both standard normal and $1$-normal.\\That ${\cal M}$ is standard normal means that the associated $\sigma$-module $(M,\phi)$ has a unique $\phi$-stable submodule of rank one on which $\phi$ acts on a basis element by multiplication with a unit in $\widehat{A}$: the {\it unit root part} $(M_{unit},\phi_{unit})$ of $(M,\phi)$. In general, $(M_{unit},\phi_{unit})$ will not be overconvergent even if $(M,\phi)$ is overconvergent. The purpose of this section is to present another construction of $\sigma$-modules departing from $(M,\phi)$ which does preserve overconvergence and allows us to recapture $(M_{unit},\phi_{unit})$ in $\Delta(\widehat{A})$, and even certain of its twists.\\

\addtocounter{satz}{1}{\bf \arabic{section}.\arabic{satz}} For $\nu\in\mathbb{Q}$ we define the $\mathbb{C}_p$-subsets $$D^{\ge \nu}:=\{x\in\mathbb{C}_p;\quad \ord_{p}(x)\ge \nu\}$$$$D^{>\nu}:=\{x\in\mathbb{C}_p;\quad \ord_{p}(x)>\nu\}.$$We use these notations also for the natural underlying rigid spaces. Let $B(\nu)_K$ be the reduced $K$-affinoid algebra consisting of power series in the free variable $V$, with coefficients in $K$, convergent on $D^{\ge\nu}$ (viewing $V$ as the standard coordinate). Thus $$B(\nu)_K=\{\sum_{\alpha\in\mathbb{N}_0}c_{\alpha}V^{\alpha};\quad c_{\alpha}\in K, \lim_{\alpha\to\infty}(\ord_{p}(c_{\alpha})+\nu\alpha)=\infty\}.$$

\addtocounter{satz}{1}{\bf \arabic{section}.\arabic{satz}}\newcounter{ri1}\newcounter{ri2}\setcounter{ri1}{\value{section}}\setcounter{ri2}{\value{satz}} Fix $\nu\in\mathbb{Q}$ and let $$B:=(B(\nu)_K)^0:=\{\sum_{\alpha\in\mathbb{N}_0}c_{\alpha}V^{\alpha}\in B(\nu)_K;\quad \ord_{p}(c_{\alpha})+\nu\alpha\ge0\mbox{ for all }\alpha\},$$ $$J:=\{q:I_1\to\mathbb{N}_0;\quad q(i)=0\mbox{ for almost all }i\in I_1\},$$$$C:=(A\widehat{\otimes}_RB)^J=\prod_JA\widehat{\otimes}_RB.$$Define a multiplication in $C$ as follows. Given $\beta=(\beta_q)_{q\in J}$ and $\beta'=(\beta'_q)_{q\in J}$ in $C$, the component at $q\in J$ of the product $\beta\beta'$ is defined as\label{multdef} $$(\beta\beta')_q=\sum_{\stackrel{(q_1,q_2)\in J^2}{q_1+q_2=q}}\beta_{q_1}\beta'_{q_2}.$$$C$ is $p$-adically complete. For $c\in\mathbb{N}_0$ we defined $[0,c]_B$ in \arabic{ovcoc1}.\arabic{ovcoc2}, and now we let$$C_c:=([0,c]_{B})^J=\prod_J[0,c]_{B},$$a complete subring of $C$. We view $C$ as a $A\widehat{\otimes}_RB$-algebra by means of the ring morphism $h:A\widehat{\otimes}_RB\to C$ defined for $y\in A\widehat{\otimes}_RB$ by $h(y)_q=y\in A\widehat{\otimes}_RB$ if $q\in J$ is the zero map $I_1\to\mathbb{N}_0$, and by $h(y)_q=0\in A\widehat{\otimes}_RB$ for all other $q\in J$. In turn, $$C\cong A\widehat{\otimes}_RB[[I_1]],$$the free power series ring on the set $I_1$ (viewed as a set of free variables).\\

\addtocounter{satz}{1}{\bf \arabic{section}.\arabic{satz}}\newcounter{einb1}\newcounter{einb2}\setcounter{einb1}{\value{section}}\setcounter{einb2}{\value{satz}} Let $\mu:S_1\to S_2$ be a homomorphism of arbitrary $R$-modules. With $I$, $i_0$, $I_1$ and $J$ from above we now define a homomorphism $$\lambda(\mu):(S_1)^I=\prod_IS_1\to (S_2)^J=\prod_JS_2.$$ Given $a=(a_i)_{i\in I}\in\prod_IS_1$, the $q$-component $\lambda(\mu)(a)_q$ of $\lambda(\mu)(a)$, for $q\in J$, is defined as follows. If $q\in J$ is the zero map $I_1\to\mathbb{N}_0$, then $\lambda(\mu)(a)_q=\mu(a_{i_0})\in S_2$. If there is a $i\in I_1$ such that $q(i)=1$ and $q(i')=0$ for all $i'\in I_1-\{i\}$, then  $\lambda(\mu)(a)_q=\mu(a_{i})\in S_2$ (for this $i$). For all other $q\in J$ we let $\lambda(\mu)(a)_q=0\in S_2$.\\Returning to the situation in \arabic{ri1}.\arabic{ri2}, the natural inclusion $\tau:\widehat{A}\to A\widehat{\otimes}_RB=\widehat{A}\widehat{\otimes}_RB$ gives us an embedding of $\widehat{A}$-modules\label{lambdadef}$$\lambda=\lambda(\tau):\widehat{A}^I=\prod_I\widehat{A}\to C=\prod_JA\widehat{\otimes}_RB.$$It is clear that $\lambda(([0,c]_R)^I)\subset C_c$.\\

\addtocounter{satz}{1}{\bf \arabic{section}.\arabic{satz}}\newcounter{lide1}\newcounter{lide2}\setcounter{lide1}{\value{section}}\setcounter{lide2}{\value{satz}} Now let ${\cal M}=(a_{i_1,i_2})_{i_1,i_2\in I}$ be a nuclear and $1$-normal $I\times I$-matrix over $\widehat{A}$. Then $$\mu:=\inf(\{\ord_p(a_{i_0,i_0}-1)\}\cup\{\ord_{p}(a_{i_1,i_2});\,(i_1,i_2)\ne(i_0,i_0)\})\ge \frac{1}{e}>0.$$If $\mu>\frac{1}{p-1}$ choose $\nu\in\mathbb{Q}$ such that $\nu>\frac{1}{p-1}-\mu$. If only $\mu\ge \frac{1}{p^b}\frac{1}{p-1}$ for some $b\in\mathbb{N}_0$ choose $\nu\in\mathbb{Q}$ such that $\nu>b$. With this $\nu$ define $B$ and $C$ as above.\\We view ${\cal M}$ as the set, indexed by $i_2\in I$, of its columns $$a_{(i_2)}:=(a_{i_1,i_2})_{i_1\in I}\in\widehat{A}^I.$$ For each $r\in\mathbb{Z}$ we now define \label{deflimm} a $J\times J$-matrix ${\cal B}^r({\cal M})=(b_{q_1,q_2}^{(r)})_{q_1,q_2\in J}$ over $A\widehat{\otimes}_RB$ associated with ${\cal M}$. To define ${\cal B}^r({\cal M})$ it is enough to define the set, indexed by $q_2\in J$, of the columns$$b_{(q_2)}^{(r)}=(b^{(r)}_{q_1,q_2})_{q_1\in J}\in \prod_JA\widehat{\otimes}_RB=C$$ of ${\cal B}^r({\cal M})$. Using the ring structure of $C$ we define$$b_{(q_2)}^{(r)}:=\lambda(a_{(i_0)})^V\lambda(a_{(i_0)})^{r}\frac{\prod_{i\in I_1}\lambda(a_{(i)})^{q_2(i)}}{\lambda(a_{(i_0)})^{|q_2|}}.$$Here $|q|=\sum_{i\in I_1}q(i)$ for $q\in J$, and $\lambda(a_{(i_0)})^V\in C$ is defined as$$\lambda(a_{(i_0)})^V:=\exp(V\log(\lambda(a_{(i_0)}))).$$For this to make sense note that $\ord_p(\lambda(a_{(i_0)})-1_C)\ge\mu>\frac{1}{p-1}$ (resp. $\ord_p(\lambda(a_{(i_0)})-1_C)\ge\mu\ge \frac{1}{p^b}\frac{1}{p-1}$), hence $\ord_p(\log(\lambda(a_{(i_0)})))\ge\mu$ (resp. $\ord_p(\log(\lambda(a_{(i_0)})))\ge \frac{1}{p-1}-b$) by \ref{logexp}(ii). Thus $V\log(\lambda(a_{(i_0)}))$ is, in view of our choice of $\nu$, indeed an element of $A\widehat{\otimes}_RB$, with $\ord_p(V\log(\lambda(a_{(i_0)})))\ge\mu+\nu>\frac{1}{p-1}$ (resp. $\ord_p(V\log(\lambda(a_{(i_0)})))\ge\frac{1}{p-1}-b+\nu>\frac{1}{p-1}$), so we can apply \ref{logexp}(i) to it.\\If the free variable $V$ specializes to integer values, $\lambda(a_{(i_0)})^V$ specializes to the usual exponentiation by integers of the unit $\lambda(a_{(i_0)})$ in $C$ (just as we use usual exponentiation for the other factors in the above definition of $b_{(q_2)}^{(r)}$). Let ${\cal B}^r_{-}({\cal M})$ be the matrix obtained from ${\cal B}^r({\cal M})$ by replacing $V$ with $-V$ (i.e. the matrix defined by the same recipe, but now using $\lambda(a_{(i_0)})^{-V}$ in place of $\lambda(a_{(i_0)})^V$ as the first factor of $b_{(q_2)}^{(r)}$).\\

\addtocounter{satz}{1}{\bf \arabic{section}.\arabic{satz}} The particular choice of $\nu$ made in \arabic{lide1}.\arabic{lide2} will play no role in the sequel. However, there is some theoretical interest in taking $\nu$ as small as possible: the smaller $\nu$, the larger $D^{\ge\nu}$ which is the parameter space for our families of $\sigma$-modules defined by the matrices ${\cal B}^r({\cal M})$ and ${\cal B}^r_{-}({\cal M})$. The ultimate result \ref{merwei} on the family of twisted unit root $L$-functions does not depend on the choice (in the prescribed range) of $\nu$ here: for \ref{merwei} it is not important how far the family extends, we only need to extend it to $D^{\ge\nu}$ for {\it some} $\nu<\infty$. But we get trace formulas, which are important for a further qualitative study, only for those members of this family of twisted unit root $L$-functions whose parameters (=locally $K$-analytic characters)  are in $D^{\ge\nu}$.\\  

\begin{pro}\label{nucov} The matrices ${\cal B}^r({\cal M})$ and ${\cal B}^r_{-}({\cal M})$ are nuclear. If ${\cal M}$ is nuclear overconvergent, then ${\cal B}^r({\cal M})$ and ${\cal B}^r_{-}({\cal M})$ are nuclear overconvergent.
\end{pro}

{\sc Proof:} Nuclearity: Given $M>0$, we need to show $\ord_{\pi}(b_{(q_2)}^{(r)})>M$ for all but finitely many $q_2\in J$. It is clear that $\ord_{\pi}(\lambda(a_{(i_0)})^V\lambda(a_{(i_0)})^m)=0$ for all $m\in\mathbb{Z}$, therefore we need to concentrate only on the factors $\prod_{i\in I_1}\lambda(a_{(i)})^{q_2(i)}$. By nuclearity of ${\cal M}$ we know that $\ord_{\pi}(\lambda(a_{(i)}))=\ord_{\pi}(a_{(i)})>M$ for all but finitely many $i\in I_1$. Therefore we need to concentrate only on those $q_2$ with support inside this finite exceptional subset of $I_1$. Among these $q_2$ we have $|q_2|>M$ for all but finitely many $q_2$. But $|q_2|>M$ (and $\ord_{\pi}(a_{(i)})\ge1$ for all $i\in I_1$) implies $$\ord_{\pi}(\prod_{i\in I_1}\lambda(a_{(i)})^{q_2(i)})=\sum_{i\in I_1}q_2(i)\ord_{\pi}(a_{(i)})\ge\sum_{i\in I_1}q_2(i)=|q_2|>M.$$Nuclearity is established. Now assume ${\cal M}$ is nuclear overconvergent. Then it can be lifted to a nuclear, overconvergent and $1$-normal matrix $\widetilde{\cal M}=(\widetilde{a}_{i_1,i_2})_{i_1,i_2\in I}$ with entries in $R[X]^{\dagger}$. Then, perhaps increasing the $c$ from our nuclearity condition, there is a $c\in\mathbb{N}$ such that $\widetilde{a}_{i_1,i_2}\in[0,c]_R$ for all $(i_1,i_2)\ne(i_0,i_0)$, and also $\widetilde{a}_{i_0,i_0}-1\in[0,c]_R$. Then all entries of ${\cal B}^r(\widetilde{\cal M})$ are in $[0,c]_{B}$. Hence ${\cal B}^r(\widetilde{\cal M})$ is nuclear and overconvergent. Clearly it is a lifting of ${\cal B}^r({\cal M})$, so we are done.\\

\addtocounter{satz}{1}{\bf \arabic{section}.\arabic{satz}}\newcounter{modlim1}\newcounter{modlim2}\setcounter{modlim1}{\value{section}}\setcounter{modlim2}{\value{satz}} Now let us look at the $\sigma$-module over $A\widehat{\otimes}_RB$ defined by the matrix ${\cal B}^r({\cal M})$. By construction, this is the $A\widehat{\otimes}_RB$-module $C$ (which in fact even is a $A\widehat{\otimes}_RB$-algebra), with the $\sigma$-linear endomorphism defined by ${\cal B}^r({\cal M})$. We view it as an analytic family, parametrized by the rigid space $\spm(B(\nu)_K)=D^{\ge\nu}$, of nuclear $\sigma$-modules over $\widehat{A}$; its fibres at points $\mathbb{Z}\cap D^{\ge\nu}$ are Wan's "limiting modules" \cite{warko}. Yet another description is due to Coleman \cite{cormk}, which we now present (in a slightly generalized form). It will be used in the proof of \ref{funclim}. \label{limmod}The nuclear matrix ${\cal M}$ over $\widehat{A}$ is the matrix in a formal basis $\{e_i\}_{i\in I}$ of a $\sigma$-linear endomorphism $\phi$ on a $\widehat{A}$-module $M$. The element $e=e_{i_0}\in M$ can also be viewed as an element of the symmetric $\widehat{A}$-algebra $\sym_{\widehat{A}}(M)$ defined by $M$, so it makes sense to adjoin its inverse to $\sym_{\widehat{A}}(M)$. Let $D$ be the subring of degree zero elements in $\sym_{\widehat{A}}(M)[\frac{1}{e}]\widehat{\otimes}_RB$: the $A\widehat{\otimes}_RB$-sub-algebra of $\sym_{\widehat{A}}(M)[\frac{1}{e}]\widehat{\otimes}_RB$ generated by all $\frac{m}{e}$ for $m\in M$. Let ${\cal I}\subset D$ be the ideal generated by all elements $\frac{m}{e}$ for $m\in M$ with $\phi(m)\in\pi M$, and let $B^r(M)$ be the $(\pi,{\cal I})$-adic completion of $D$. For all $\alpha\in (A\widehat{\otimes}_RB)^{\times}$, all $m_1,m_2\in M$, if we set $e'=\alpha e+\pi m_1$, we have \begin{gather}\frac{m_2}{e'}=\frac{m_2}{\alpha e}\sum_{i=0}^{\infty}(\frac{\pi}{\alpha}\frac{m_1}{e})^i\tag{$*$}\end{gather}in $B^r(M)$. By our assumptions on ${\cal M}$ we know $\phi(e)-e\in\pi M$. Therefore there exists a unique $\sigma$-linear ring endomorphism $\psi$ of $B^r(M)$ with $\psi(\frac{m}{e})=\frac{\phi(m)}{\phi(e)}$ for all $m\in M$: Take $(*)$ as a definition, with $e'=\phi(e)$, $m_2=\phi(m)$ and $\alpha=1$. Similarly as in \arabic{lide1}.\arabic{lide2} we can define, for integers $r\in\mathbb{Z}$, the element $$(\frac{\phi(e)}{e})^{V+r}=\exp(V\log(\frac{\phi(e)}{e}))(\frac{\phi(e)}{e})^r$$ of $B^r(M)$. We define the $\sigma$-linear endomorphism $B^r(\phi)$ of $B^r(M)$ by$$B^r(\phi)(y)=(\frac{\phi(e)}{e})^{V+r}\psi(y)$$for $y\in B^r(M)$. Clearly ${\cal B}^r({\cal M})$ is the matrix of $B^r(\phi)$ acting on the formal basis$$\{\prod_{i\in I}(\frac{e_i}{e})^{q(i)}\}_{q\in J}$$of $B^r(M)$ over $A\widehat{\otimes}_RB$. The $\sigma$-module defined by ${\cal B}^r_{-}({\cal M})$ is described similarly.\\ 

\begin{pro}\label{funclim} The $\sigma$-similarity classes (over $A\widehat{\otimes}_RB$) of ${\cal B}^r({\cal M})$ and ${\cal B}^r_{-}({\cal M})$ depend only on the $\sigma$-similarity class (over $\widehat{A}$) of ${\cal M}$.
\end{pro}

{\sc Proof:} We prove this for ${\cal B}^r({\cal M})$, the argument for ${\cal B}^r_{-}({\cal M})$ is the same. It is enough to prove that $B^r(M)$, as a $A\widehat{\otimes}_RB$-module together with its $\sigma$-linear endomorphism $B^r(\phi)$, depends only on the $\sigma$-module $M$. Let ${\cal M}'=(a'_{i_1,i_2})_{{i_1,i_2}\in I}$ be another $1$-normal nuclear matrix over $\widehat{A}$ which is $\sigma$-similar to ${\cal M}$. We can view ${\cal M}'$ as the matrix of the {\it same} $\sigma$-linear endomorphism $\phi$ on the {\it same} $\widehat{A}$-module $M$, but in another formal basis $\{e'_i\}_{i\in I}$. For the element $e'=e'_{i_0}$ our assumptions imply $\phi(e')-e'\in\pi M$. Therefore $e'$ and $e$ both generate the unit root part modulo $\pi$ of $M$, hence there is a $\alpha\in\widehat{A}^{\times}$ with $e'-\alpha e\in\pi M$. Observe that$$\alpha e\equiv e'\equiv\phi(e')\equiv\phi(\alpha e)\equiv\sigma(a)\phi(e)\equiv\sigma(\alpha)e$$modulo $\pi M$. Since $R^{\times}$ is the subgroup of $\widehat{A}^{\times}$ fixed by $\sigma$ we may and will assume $\alpha\in{R}^{\times}$. From $(*)$ in \arabic{modlim1}.\arabic{modlim2} it follows that for $m\in M$ the element $\frac{m}{e'}$ of the $\pi$-adic completion of $\sym_{\widehat{A}}(M)[M^{-1}]\widehat{\otimes}_RB$ actually lies in its subring $B^r(M)$. By a symmetry argument we deduce that $B^r(M)$ is the same when constructed with respect to $e$ or with respect to $e'$. Moreover the endomorphism $\psi$ on $B^r(M)$ is the same when constructed with respect to $e$ or with respect to $e'$: it is uniquely determined by its action on $B^r(M)\cap \sym_{\widehat{A}}(M)[M^{-1}]\widehat{\otimes}_RB$, where it is characterized by $\psi(\frac{m_1}{m_2})=\frac{\phi(m_1)}{\phi(m_2)}$ for $m_1, m_2\in M$. Now let$$B^r(\phi)'(y)=(\frac{\phi(e')}{e'})^{V+r}\psi(y)$$for $y\in B^r(M)$. The needed $A\widehat{\otimes}_RB$-linear endomorphism $\lambda_r$ of $B^r(M)$ satisfying $\lambda_r\circ B^r(\phi)'=B^r(\phi)\circ\lambda_r$ we now define to be the multiplication with $(\frac{e'}{\alpha e})^{V+r}\in B^r(M)$ (by now obviously defined). Here we use that $\alpha\in{R}^{\times}$.\\

Now suppose ${\cal M}=(a_{i_1,i_2})_{i_1,i_2\in I}$ is even a standard $1$-normal nuclear $\widehat{A}$-matrix. Define ${\cal M}_{unit}:=a_{i_0,i_0}\in\widehat{A}$ and $({\cal M}_{unit})^V=\exp(V\log({\cal M}_{unit}))\in A\widehat{\otimes}_RB$ as in \arabic{lide1}.\arabic{lide2}.\\

\begin{satz}\label{rk1res} For $s\in\mathbb{Z}$ we have the following equalities in $\Delta(A\widehat{\otimes}_RB)$:$$[({\cal M}_{unit})^s({\cal M}_{unit})^V]=\bigoplus_{r\ge1}(-1)^{r-1}r[{\cal B}^{s-r}({\cal M})\otimes\bigwedge^r({\cal M})]$$$$[({\cal M}_{unit})^s({\cal M}_{unit})^{-V}]=\bigoplus_{r\ge1}(-1)^{r-1}r[{\cal B}^{s-r}_{-}({\cal M})\otimes\bigwedge^r({\cal M})].$$
\end{satz}

{\sc Proof:} We prove the first equality, the second is proved similarly. First note that our assumptions imply that $\pi^{r-1}$ divides $\bigwedge^r({\cal M})$, so the right hand side converges. Since ${\cal M}$ is standard $1$-normal we have ${\cal B}^{s-r}({\cal M})=({\cal M}_{unit})^s{\cal B}^{-r}({\cal M})$ so we may assume $s=0$. Let ${\cal M}''=(a''_{i_1,i_2})_{i_1,i_2\in I}$ be the $\widehat{A}$-matrix with $a''_{i_1,i_2}=a_{i_1,i_2}$ for all $(i_1,i_2)\in (I_1\times I_1)\cup\{(i_0,i_0)\}$, and $a''_{i_1,i_2}=0$ for the other $(i_1,i_2)$. Since ${\cal M}$ is standard $1$-normal we see that $[{\cal B}^{-r}({\cal M})\otimes\bigwedge^r({\cal M})]=[{\cal B}^{-r}({\cal M}'')\otimes\bigwedge^r({\cal M}'')]$ in view of the relations divided out in the definition of $\Delta(A\widehat{\otimes}_RB)$. Hence we may assume ${\cal M}={\cal M}''$. Suppose that $i_0$ is minimal in the ordering of $I$ (which we tacitly chose to define $\bigwedge^r(I)$ and $\bigwedge^r({\cal M})$, see \arabic{matop1}.\arabic{matop2}). For $r\ge 1$ let $M_r$ be the $A\widehat{\otimes}_RB$-module $(A\widehat{\otimes}_RB)^{(J\times\bigwedge^r(I))}$. It has the formal basis $(e_{(q,\vec{\imath})})_{q\in J, \vec{\imath}\in\bigwedge^r(I)}$, where $e_{(q,\vec{\imath})}$ is the $(q,\vec{\imath})$-th column of the identity $(J\times\bigwedge^r(I))\times(J\times\bigwedge^r(I))$-matrix. The matrix ${\cal B}^{-r}({\cal M})\otimes\bigwedge^r({\cal M})$ describes the action of a $\sigma$-linear endomorphism $\phi_r$ of $M_r$ on this basis. Actually we will need $r$ copies of $(M_r,\phi_r)$ and its formal basis $(e_{(q,\vec{\imath})})_{q\in J, \vec{\imath}\in\bigwedge^r(I)}$: We denote them by $(M_r^{(\ell)},\phi^{(\ell)}_r)$ and $(e^{(\ell)}_{(q,\vec{\imath})})_{(q,\vec{\imath})}$ for $1\le\ell\le r$. We get the $\sigma$-module $$(M^{\bullet}_r,\phi_r^{\bullet})=\oplus_{1\le\ell\le r}(M_r^{(\ell)},\phi^{(\ell)}_r)$$with formal basis$$H_r=(e^{(\ell)}_{(q,\vec{\imath})})_{q\in J, \vec{\imath}\in\bigwedge^r(I),\ell\in\{1,\ldots,r\}}.$$Define $A\widehat{\otimes}_RB$-linear maps$$\alpha^{(\ell)}_r:M_r^{(\ell)}\to M_{r+1}^{(\ell)}$$$$\beta^{(\ell)}_r:M_r^{(\ell)}\to M_{r+1}^{(\ell+1)}$$as follows. For $\vec{\imath}=(i_1,\ldots,i_r)\in\bigwedge^r(I)$ with $i_1<\ldots<i_r$, and another $i\in I$, let $\tau(\vec{\imath},i)=\max(\{t\le r;i_t<i\}\cup\{0\})$, and if in addition $i\ne i_{\tau(\vec{\imath},i)+1}$ let$$[\vec{\imath},i]=(i_1,\ldots,i_{\tau(\vec{\imath},i)},i,i_{\tau(\vec{\imath},i)+1},\ldots,i_r)\in\bigwedge^{r+1}(I).$$For $q\in J$ and $i\in I_1$ with $q(i)\ne0$ define $q^{i-}\in J$ by $q^{i-}(i')=q(i')$ for $i'\in I_1-\{i\}$, and $q^{i-}(i)=q(i)-1$. Now set$$\alpha^{(\ell)}_r(e^{(\ell)}_{(q,\vec{\imath})})=e^{(\ell)}_{(q,[\vec{\imath},i_0])}$$if $i_1\ne i_0$, and set $\alpha^{(\ell)}_r(e^{(\ell)}_{(q,\vec{\imath})})=0$ if $i_1=i_0$. Set $$\beta^{(\ell)}_r(e^{(\ell)}_{(q,\vec{\imath})})=\sum_t(-1)^{\tau(\vec{\imath},i_t)}e^{(\ell+1)}_{(q^{i_t-},[\vec{\imath},i_t])}$$where the sum runs through all $1\le t\le r$ with $i_t\ne i_0$, with $i_t\ne i_{\tau(\vec{\imath},i_t)+1}$ and with $q(i_t)\ne0$. One checks that $\phi^{(\ell)}_{r+1}\circ\alpha^{(\ell)}_r=\alpha^{(\ell)}_r\circ\phi^{(\ell)}_r$ (use the standard $1$-normality of ${\cal M}$), and that $\phi^{(\ell+1)}_{r+1}\circ\beta^{(\ell)}_r=\beta^{(\ell)}_r\circ\phi^{(\ell)}_r$ (use ${\cal M}={\cal M}''$). Hence for $$\psi^{\bullet}_r=\oplus_{1\le\ell\le r}(\alpha^{(\ell)}_r\oplus\beta^{(\ell)}_r):M^{\bullet}_r\to M^{\bullet}_{r+1}$$we have $\phi^{\bullet}_{r+1}\circ\psi^{\bullet}_r=\psi^{\bullet}_r\circ\phi^{\bullet}_r$. Also note that $\phi_1^{\bullet}=\phi_1^{(1)}$ on $M_1^{\bullet}=M_1^{(1)}$ restricts on the rank one $A\widehat{\otimes}_RB$-submodule $M_0^{\bullet}$ spanned by the basis element $e_{(0,i_0)}^{(1)}\in J\times I=J\times\bigwedge^1(I)$ to a $\sigma$-linear endomorphism $\phi_0$ with matrix $({\cal M}_{unit})^V$. Let $\psi_0^{\bullet}:M_0^{\bullet}\to M_1^{\bullet}$ be the inclusion and consider\begin{gather}0\to M_0^{\bullet}\stackrel{\psi_0^{\bullet}}{\to}M_1^{\bullet}\stackrel{\psi_1^{\bullet}}{\to}M_2^{\bullet}\stackrel{\psi_2^{\bullet}}{\to}\ldots.\tag{$*$}\end{gather}We saw that this sequence is equivariant for the $\sigma$-linear endomorphisms $\phi_r^{\bullet}$ which are described by matrices as occur in the statement of the theorem, so it remains to show that $(*)$ is split exact; more precisely, that for each $r$ there are disjoint subsets $G_r^1$ and $G_r^2$ of $M_r^{\bullet}$ with the following properties: $\psi_r^{\bullet}$ induces a bijection of sets $G_r^2\cong G_{r+1}^1$, and the union $G_r^1\cup G_r^2$ is a formal basis for $M_r^{\bullet}$ (transforming under an invertible matrix to the formal basis $H_r$). We let $G_0^1=\emptyset$, $G_0^2=H_0=\{e^{(1)}_{(0,i_0)}\}$. For $r\ge1$ we let$$G^1_r=\{\psi^{\bullet}_{r-1}(h); h\in H_{r-1}\}.$$We let $G_r^2$ be the subset of $H_r$ consisting of those $e^{(\ell)}_{(q,\vec{\imath})}$ with $\ell, q$ and $\vec{\imath}=(i_1,\ldots,i_r)\in\bigwedge^r(I)$ satisfying one of the following conditions: either$$[\ell=1\mbox{ and }((i_1\ne i_0)\mbox{ or }(i_1=i_0\mbox{ and }\exists(i\in I_1-\{i_2,\ldots,i_r\}):\quad q(i)\ne0))]$$or$$[\ell\ne 1\mbox{ and }((i_1\ne i_0\mbox{ and }\forall(1\le k\le r)\exists(i\in I_1-\{i_k\}):\quad q(i)\ne0)$$$$\mbox{ or }(i_1=i_0\mbox{ and }\exists(i\in I_1-\{i_2,\ldots,i_r\}):\quad q(i)\ne0))].$$The desired properties are formally verified, the proof is complete.\\

\begin{kor}\label{wafa} Suppose our ${\cal M}$ is also overconvergent nuclear. Then for each $s\in\mathbb{Z}$ the series$$\prod_{\overline{x}\in X}\frac{1}{\det(1-({\cal M}_{unit})^s_{\overline{x}}({\cal M}_{unit})_{\overline{x}}^yT^{\deg(\overline{x})})}$$defines a meromorphic function in the variables $T$ and $y$ on $\mathbb{A}_{\mathbb{C}_p}^1\times D^{\ge\nu}$, specializing for $y\in D^{\ge\nu}(K)$ to $L({\cal M}_{unit}^{s+y},T)$.
\end{kor}

{\sc Proof:} The series is trivially holomorphic on $D^{>0}\times D^{\ge\nu}$. We claim that it is equal to$$\prod_{r\ge1}(\prod_{i=0}^d\det(1-\psi[{\cal B}^{s-r}({\cal M})\otimes\bigwedge^r({\cal M})\otimes{\mathcal D}^{\wedge i}]T)^{(-1)^{i-1}})^{(-1)^{r-1}r}$$which clearly extends as desired. It suffices to prove equality at all specializations $V=y$ at $K$-rational points $y\in D^{\ge\nu}(K)$ (since these $y$ are Zariski dense in $D^{\ge\nu}$). But for such $y$ both series coincide with
$$\prod_{r\ge1}L({\cal B}^{s-r}({\cal M})|_{V=y}\otimes\bigwedge^r({\cal M}),T)^{(-1)^{r-1}r}:$$For the series in the statement of \ref{wafa} this follows from \ref{lfndbkt} and \ref{rk1res}, for the first series written in this proof this follows from \ref{tracf}. 

\section{Weight space ${\cal W}$}

In this section we describe a $K$-rigid analytic space ${\cal W}$ whose set of $\mathbb{C}_p$-valued points can be identified with the set of locally $K$-analytic characters $\kappa:R^{\times}\to\mathbb{C}_p$ occuring in Theorem \ref{hauptsatz}.\\ 

\addtocounter{satz}{1}{\bf \arabic{section}.\arabic{satz}} For a $K$-analytic group manifold $G$ (see \cite{bou}) we denote by $\ho_{K\mbox{-}an}(G,\mathbb{C}_p^{\times})$ the group of locally $K$-analytic characters $G\to\mathbb{C}_p^{\times}$: characters which locally on $G$ can be expanded into power series in $\dim_K(G)$-many variables. If $K=\mathbb{Q}_p$ these are precisely the {\it continuous} characters $G\to\mathbb{C}_p^{\times}$. The examples relevant for us are $G=R$, $G=R^{\times}$, $G=U_{R}^{(1)}$ and $G=\overline{U}_{R}^{(1)}$, where we write$$U_{R}^{(1)}:=1+\pi R\quad\mbox{ and }\quad \overline{U}_{R}^{(1)}:=\frac{U_{R}^{(1)}}{(U_{R}^{(1)})_{\mbox{tors}}}.$$To extinguish any confusion, although in these examples $G$ even carries a natural structure of $K$-{\it rigid} group variety, the definition of $\ho_{K\mbox{-}an}(G,\mathbb{C}_p^{\times})$ does not refer to this (indeed more "rigid") structure: local $K$-analycity of a character $\kappa$ requires only that $\kappa$, as a $\mathbb{C}_p^{\times}\subset\mathbb{C}_p$ valued function on $G$, can be expanded into convergent power series on each member of some open covering of $G$ --- an open covering in the naive sense, not necessarily admissible in the sense of rigid geometry.\\

\addtocounter{satz}{1}{\bf \arabic{section}.\arabic{satz}}\newcounter{weicons1}\newcounter{weicons2}\setcounter{weicons1}{\value{section}}\setcounter{weicons2}{\value{satz}} Let ${\cal G}={\cal G}_{\pi}$ be the Lubin-Tate formal group over $R$ corresponding to our chosen uniformizer $\pi\in R$ (see \cite{lan}). For $x\in R$ denote by $[x]\in U.R[[U]]$ the formal power series which defines the multiplication with $x$ in the formal $R$-module ${\cal G}$. The $R$-module $\ho_{{\cal O}_{\mathbb{C}_p}}({\cal G}\widehat{\otimes}{{\cal O}_{\mathbb{C}_p}},\mathbb{G}_{m,{\cal O}_{\mathbb{C}_p}})$ is free of rank one. Fix a generator with corresponding power series $F(Z)\in Z.{\cal O}_{\mathbb{C}_p}[[Z]]$. Substitution yields power series $F([x])\in U.{\cal O}_{\mathbb{C}_p}[[U]]$ for $x\in R$. By \cite{scte} we have a group isomorphism$$D^{>0}\stackrel{\cong}{\longrightarrow}\ho_{K\mbox{-}an}(R,\mathbb{C}_p^{\times})$$$$z\mapsto[x\mapsto 1+F([x])(z)].$$Here $D^{>0}$ carries the group structure defined by ${\cal G}$. Let $m\in\mathbb{Z}_{\ge-1}$ be minimal such that $\pi^m\log(U_{R}^{(1)})\subset R$. Since $(U^{(1)}_R)_{\mbox{tors}}=\ke(\log)$ we have a well defined injective homomorphism of $K$-analytic group varieties$$\overline{U}_{R}^{(1)}\stackrel{\theta}{\longrightarrow}R,\quad u\mapsto\pi^m\log(u)=\theta(u)$$inducing a homomorphism$$\ho_{K\mbox{-}an}(R,\mathbb{C}_p^{\times})\stackrel{\delta}{\longrightarrow}\ho_{K\mbox{-}an}(\overline{U}_{R}^{(1)},\mathbb{C}_p^{\times}).$$Note that $\koke(\theta)$ is finite. Thus $\ke(\delta)$ is finite, and on the other hand $\delta$ is surjective (since $\mathbb{C}_p^{\times}$ is divisible). In other words, $\ho_{K\mbox{-}an}(\overline{U}_{R}^{(1)},\mathbb{C}_p^{\times})$ is the quotient of $\ho_{K\mbox{-}an}(R,\mathbb{C}_p^{\times})$ by a finite subgroup $\Delta\subset\ho_{K\mbox{-}an}(R,\mathbb{C}_p^{\times})$. The formal group law ${\cal G}$ defines a structure of $\mathbb{C}_p$-rigid analytic group variety on $D^{>0}$ (with its standard coordinate $U$). By means of the above isomorphism we view $\ho_{K\mbox{-}an}(R,\mathbb{C}_p^{\times})$ as its group of $\mathbb{C}_p$-valued points. Accordingly, we view $\ho_{K\mbox{-}an}(\overline{U}_{R}^{(1)},\mathbb{C}_p^{\times})$ as the group $(D^{>0}/\Delta)(\mathbb{C}_p)$ of $\mathbb{C}_p$-valued points of the $\mathbb{C}_p$-rigid group variety $D^{>0}/\Delta$. Let $$U_{R}^{(1)}\stackrel{u}{\leftarrow}R^{\times}\stackrel{v}{\to}\mathbb{\mu}_{q-1}$$be the natural projections. We have $(U_{R}^{(1)})_{\mbox{tors}}=\mathbb{\mu}_{p^a}$ for some $a\ge 0$. Let $${\cal F}=\{0,\ldots,p^a-1\}\times\{0,\ldots,q-2\}.$$For each $(s,t)\in{\cal F}$ let ${\cal W}_{(s,t)}$ be a copy of the $\mathbb{C}_p$-rigid space $D^{>0}/\Delta$, and let $${\cal W}:=\coprod_{(s,t)\in{\cal F}}{\cal W}_{(s,t)}.$$For an element $\omega\in{\cal W}_{(s,t)}(\mathbb{C}_p)\subset {\cal W}(\mathbb{C}_p)$ we define the character$$\kappa_{\omega}:R^{\times}\to \mathbb{C}_p,\quad r\mapsto u(r)^s\widetilde{\omega}(u(r))v(r)^t=:\kappa_{\omega}(r)$$where $\widetilde{\omega}\in\ho_{K\mbox{-}an}({U}_{R}^{(1)},\mathbb{C}_p^{\times})$ is the image of $\omega$ under the natural map$${\cal W}_{(s,t)}(\mathbb{C}_p)\cong (D^{>0}/\Delta)(\mathbb{C}_p)\cong \ho_{K\mbox{-}an}(\overline{U}_{R}^{(1)},\mathbb{C}_p^{\times})\to \ho_{K\mbox{-}an}({U}_{R}^{(1)},\mathbb{C}_p^{\times}).$$Since $R^{\times}=\mathbb{\mu}_{q-1}\times U^{(1)}_R$ we get:

\begin{pro}\label{weigsp} The assignment $\omega\mapsto\kappa_{\omega}$ defines a bijection$${\cal W}(\mathbb{C}_p)\cong\ho_{K\mbox{-}an}(R^{\times},\mathbb{C}_p^{\times}).$$Thus $\ho_{K\mbox{-}an}(R^{\times},\mathbb{C}_p^{\times})$ can be viewed as the set of $\mathbb{C}_p$-valued points of the $\mathbb{C}_p$-rigid variety ${\cal W}$.
\end{pro}

\begin{lem}\label{lubtat} For $\nu\in\mathbb{Q}$, $\nu>\frac{m}{e}+\frac{1}{p-1}$, there exists an open embedding of $\mathbb{C}_p$-rigid varieties $\iota:D^{\ge\nu}\to D^{>0}$ such that for all $x\in R$ and all $y\in D^{\ge\nu}$ we have $$1+F([x])(\iota(y))=\exp(\pi^{-m}xy).$$
\end{lem}

{\sc Proof:} Let $\log_{\cal G}$ be the logarithm of ${\cal G}$. Write $F(Z)=\Omega.Z+\ldots\in Z.{\cal O}_{\mathbb{C}_p}[[Z]].$ Then we have the identity of formal power series (cf. \cite{scte} sect.4)$$1+F([x])=\exp(\Omega\log_{\cal G}([x]))$$in ${\cal O}_{\mathbb{C}_p}[[U]]$. But $\log_{\cal G}([x])=x.\log_{\cal G}(U)$ by \cite{lan} 8.6 Lemma 2, therefore it is enough to find $\iota$ with$$\log_{\cal G}(\iota(y))=\pi^{-m}\Omega^{-1}y.$$By \cite{lan} 8.6 Lemma 4 the power series inverse to $\log_{\cal G}$ defines an open embedding $\exp_{\cal G}:D^{\ge\beta}\to D^{>0}$ for $\beta> \frac{1}{e(q-1)}$. Thus $\iota(y)=\exp_{\cal G}(\pi^{-m}\Omega^{-1} y)$ is appropriate; it is well defined on $D^{\ge\nu}$ because we have $\ord_{p}(\Omega)=\frac{1}{p-1}-\frac{1}{e(q-1)}$ by \cite{scte}, hence $\nu-\frac{m}{e}-\ord_{p}(\Omega)>\frac{1}{e(q-1)}$.\\

\addtocounter{satz}{1}{\bf \arabic{section}.\arabic{satz}} Example: Consider the case $K=\mathbb{Q}_p$, $\pi=p$. Then $${\cal G}=\mathbb{G}_m,\quad\quad\log_{\cal G}(Z)=\log(1+Z),\quad\quad m=-1$$$$[x]=(1+U)^x-1=\sum_{n\ge1}{x \choose n }U^n\quad\quad(x\in\mathbb{Z}_p).$$We may choose $F(Z)=Z$, and for $\nu>\frac{2-p}{p-1}$ the associated embedding$$\iota:D^{\ge\nu}\to D^{>0};\quad y\mapsto \iota(y)=\exp(py)-1$$is an isomorphism $\iota:D^{\ge\nu}\cong D^{\ge\nu+1}\subset D^{>0}$.\\

\section{Meromorphic continuation of unit root $L$-functions}
\label{mostnovel}
In this section we prove (the infinite rank version of) Theorem \ref{hauptsatz}. Let us give a sketch. For simplicity suppose that $\alpha\in\widehat{A}$ is a matrix of the ordinary unit root part of some nuclear overconvergent $\sigma$-module $M$ over $A$ (in the general case, $\alpha$ splits into two factors each of which is of this more special type and can "essentially" be treated separately). An appropriate multiplicative decomposition of $\alpha$ (see \ref{discrdec}) allows us to assume that $\alpha$ is a $1$-unit. Then the results of section \ref{resol}, together with the trace formula \ref{tracf} already show meromorphy of $L_{\alpha}$ on $\mathbb{A}\times {\cal W}^0$ for some open subspace ${\cal W}^0\subset {\cal W}$ meeting each component of ${\cal W}$: this is essentially what we proved in \ref{wafa}. More precisely we get a decomposition of $L_{\alpha}$ into holomorphic functions on $\mathbb{A}\times {\cal W}^0$ which are Fredholm determinants $\det(\psi)$ of certain completely continuous operators $\psi$ arising from limiting modules. We express the coefficients of the logarithms of these $\det(\psi)$ through the traces $\tr(\psi^f)$ of iterates $\psi^f$ of these $\psi$. Then we repeat the limiting module construction in each fibre $\overline{x}\in X$ and prove its commutation with its global counterpart. Together with the trace formula \ref{tracf} and the description of the embedding ${\cal W}^0\to{\cal W}$ given in \ref{lubtat} this can be used to show that all the functions $\tr(\psi^f)$, a priori living on ${\cal W}^0$, extend to functions on ${\cal W}$, bounded by $1$. By the general principle \ref{convext} below this implies the theorem.\\  

\begin{lem}\label{convext} For $m\in\mathbb{N}$ let $g_m(U)\in {\cal O}_{\mathbb{C}_p}[[U_1,\ldots,U_g]]$. Suppose there exists a $\tau>0$ such that$$f(T,U)=\exp(-\sum_{m=1}^{\infty}\frac{g_m(U)}{m}T^m)\in \mathbb{C}_p[[T,U_1,\ldots,U_g]]$$converges on $\mathbb{A}_{\mathbb{C}_p}^1\times(D^{\ge \tau})^g$, where $T$ resp. $U_1,\ldots,U_g$ are the standard coordinates on $\mathbb{A}_{\mathbb{C}_p}^1$, resp. on $(D^{\ge \tau})^g$. Then $f(T,U)$ converges on all of $\mathbb{A}_{\mathbb{C}_p}^1\times(D^{>0})^g$.\end{lem}

{\sc Proof:} We reduce the convergence of $f$ at a given point $x\in\mathbb{A}_{\mathbb{C}_p}^1\times(D^{>0})^g$ to the convergence of $f$ at regions -- chosen in dependence on $x$ -- of $\mathbb{A}_{\mathbb{C}_p}^1\times(D^{\ge \tau})^g$ with possibly much larger $T$-coordinates than the $T$-coordinate of $x$. For $m\ge1$ let $$I_m=\{i=(i_1,\ldots,i_m)\in(\mathbb{N}_0)^m;\quad i_1+2i_2+\ldots+mi_m=m\}.$$We may write$$f(T,U)=1+\sum_{m=1}^{\infty}\alpha_m(U)T^m$$$$\alpha_m(U)=\sum_{i\in I_m}(-1)^{i_1+\ldots+i_m}\frac{g_1(U)^{i_1}\ldots g_m(U)^{i_m}}{i_1!\ldots i_m!1^{i_1}2^{i_2}\ldots m^{i_m}}=\sum_{\ell\in(\mathbb{N}_0)^g}\beta_{m,\ell}U^{\ell}$$$$\beta_{m,\ell}=\sum_{i\in I_m}(-1)^{i_1+\ldots+i_m}\frac{\gamma_{m,\ell,(i_1,\ldots,i_m)}}{i_1!\ldots i_m!1^{i_1}2^{i_2}\ldots m^{i_m}}$$for certain $\gamma_{m,\ell,(i_1,\ldots,i_m)}\in {\cal O}_{\mathbb{C}_p}$. We have the estimate\begin{align}\ord_p(\beta_{m,l})&\ge\min_{i\in I_m}-\ord_p(i_1!\ldots i_m!1^{i_1}2^{i_2}\ldots m^{i_m})\notag \\
{} &  \ge\min_{i\in I_m}-\sum_{j=1}^m(\ord_p(i_j!)+i_j(\frac{j}{p}))\notag \\
{} &  \ge\min_{i\in I_m}-\sum_{j=1}^m(i_j+i_j(\frac{j}{p}))\notag \\
{} &  \ge\min_{i\in I_m}-2\sum_{j=1}^mji_j=-2m.\notag\end{align}Now let $(t,u_1,\ldots,u_g)\in\mathbb{A}_{\mathbb{C}_p}^1\times(D^{>0})^g$ be given. Set $$0<\rho=\min\{1,\frac{\ord_p(u_1)}{\tau},\ldots,\frac{\ord_p(u_g)}{\tau}\}\le1$$$$\lambda=\frac{\ord_p(t)-2(1-\rho)}{\rho}.$$Then we find\begin{align}\ord_p(\beta_{m,\ell}u^{\ell}t^m)&\ge-2m(1-\rho)+\rho\ord_p(\beta_{m,\ell})+\rho|\ell|\tau+\rho m\lambda+2m(1-\rho)\notag \\
{} & =\rho(\ord_p(\beta_{m,l})+|\ell|\tau+m\lambda)\notag \end{align}and this term tends to infinity as $|\ell|+m$ tends to infinity since by hypothesis $f$ converges at the points $(\widetilde{t},\widetilde{u}_1,\ldots,\widetilde{u}_g)$ with $\ord_p(\widetilde{t})\ge\lambda$ and $\ord_p(\widetilde{u}_i)\ge \tau$. The lemma follows.\\

Now let $I$ and $i_0\in I$ be as in \arabic{nordef1}.\arabic{nordef2}. In particular we can talk about $1$-normality and standard normality of $I\times I$-matrices.\\ 

\begin{lem}\label{hndec} Suppose the nuclear $I\times I$-matrix ${\cal M}$ with entries in $\widehat{A}$ is $1$-normal. Then ${\cal M}$ is $\sigma$-similar to a standard $1$-normal nuclear $I\times I$-matrix.
\end{lem}

{\sc Proof:} In case $A=R[X]^{\dagger}$, this is the translation of \cite{warko} Lemma 6.5 into matrix terminology. But the proof works for general $A$.\\

\begin{lem}\label{discrdec} Let ${\cal N}$ be a nuclear overconvergent $I\times I$-matrix over $A$ which is $\sigma$-similar to a standard normal nuclear $I\times I$-matrix over $\widehat{A}$. Then there exist a $\xi\in A$ and a nuclear overconvergent $1$-normal $I\times I$ matrix ${\cal M}$ over $A$, both unique up to $\sigma$-similarity, such that\\(i) the $1\times1$-matrix $\xi^{q-1}$ is $\sigma$-similar to $1\in A$, and\\(ii) $\xi{\cal M}$ is $\sigma$-similar to ${\cal N}$.
\end{lem}

{\sc Proof:} For the existence see Wan \cite{wahrk} (there $I$ is finite, but at this point this is not important). For the uniqueness (which by the way we do not need in the sequel) we follow Coleman \cite{cormk}. Let $\xi'$ and ${\cal M}'$ be another such pair. Then $\xi'=a\xi$ for some $a\in A^{\times}$, hence $a^{q-1}=\frac{\sigma(b)}{b}$ for some $b\in A^{\times}$ by hypothesis (i) for $\xi$ and $\xi'$. On the other hand, from hypothesis (ii) for ${\cal M}$ and ${\cal M}'$ it follows that ${\cal M}'$ and $\frac{1}{a}{\cal M}$ are $\sigma$-similar, and by $1$-normality of ${\cal M}$ and ${\cal M}'$ this implies $a=\frac{d\sigma(c)}{c}$ for some $c,d\in A^{\times}$ with $d-1\in\pi A$. Thus for $e=\frac{b}{c^{q-1}}$ we have $d^{q-1}=\frac{\sigma(e)}{e}$. In particular $(\sigma(e)-e)\in\pi A$, hence $e\in R+\pi A$, so we may assume in addition $e-1\in\pi A$. For (the unique) $f\in A$ with $f^{q-1}=e$ and $f-1\in\pi A$ we then see $d=\frac{\sigma(f)}{f}$. Thus $a=\frac{\sigma(ef)}{ef}$ and it follows that $\xi$ is $\sigma$-similar to $\sigma'$, and ${\cal M}$ to ${\cal M'}$. We are done.\\

\addtocounter{satz}{1}{\bf \arabic{section}.\arabic{satz}}\newcounter{faslim1}\newcounter{faslim2}\setcounter{faslim1}{\value{section}}\setcounter{faslim2}{\value{satz}} Let ${\cal M}$ be a standard $1$-normal nuclear $I\times I$-matrix over $\widehat{A}$. Define $I_1=I-\{i_0\}$ and $J$ as in \arabic{ri1}.\arabic{ri2}. Let $\overline{x}$ be a closed point of $X$ of degree $f$ and write ${\cal M}_{\overline{x}}=(a^{\overline{x}}_{i_1,i_2})_{i_1,i_2\in I}$ for the fibre matrix ${\cal M}_{\overline{x}}$ with entries $a^{\overline{x}}_{i_1,i_2}$ in $R_f$ as defined in \arabic{lfnd1}.\arabic{lfnd2}. We denote its $i_2$-column for $i_2\in I$ by $$a^{\overline{x}}_{(i_2)}:=(a^{\overline{x}}_{i_1,i_2})_{i_1\in I}\in\prod_IR_f.$$ Let $$\eta=1+F([\pi^m\log(a^{\overline{x}}_{i_0,i_0})])\in{\cal O}_{\mathbb{C}_p}[[U]]$$ with $F$ and $m$ as in \arabic{weicons1}.\arabic{weicons2}. For $r\in\mathbb{Z}$ we now define a nuclear $J\times J$-matrix $\widetilde{\cal B}^r({\cal M}_{\overline{x}})=(b^{(r),\overline{x}}_{q_1,q_2})_{q_1,q_2\in J}$ with entries in ${\cal O}_{\mathbb{C}_p}[[U]]$. It is enough to give the columns $$b^{(r),\overline{x}}_{(q_2)}:=(b^{(r),\overline{x}}_{q_1,q_2})_{q_1\in J}\in\prod_{j\in J}{\cal O}_{\mathbb{C}_p}[[U]],$$indexed by $q_2\in J$, of $\widetilde{\cal B}^r({\cal M}_{\overline{x}})$. The natural embedding $\rho:R_f\to {\cal O}_{\mathbb{C}_p}[[U]]$ defines a map$$\lambda=\lambda(\rho):\prod_IR_f\to\prod_J{\cal O}_{\mathbb{C}_p}[[U]]$$ as explained in \arabic{einb1}.\arabic{einb2}. We will also need the ${\cal O}_{\mathbb{C}_p}[[U]]$-algebra structure on $\prod_J{\cal O}_{\mathbb{C}_p}[[U]]$ analogous to that on $C$ in \arabic{ri1}.\arabic{ri2}. Namely, the one we get from the natural identification $$\prod_J{\cal O}_{\mathbb{C}_p}[[U]]\cong {\cal O}_{\mathbb{C}_p}[[U]][[I_1]],$$the formal power series ring over ${\cal O}_{\mathbb{C}_p}[[U]]$ on the set $I_1$ (viewed as a set of free variables). Using this ${\cal O}_{\mathbb{C}_p}[[U]]$-algebra structure we set$$b^{(r),\overline{x}}_{(q_2)}:=\eta\lambda(a_{(i_0)}^{\overline{x}})^{r}\frac{\prod_{i\in I_1}\lambda(a_{(i)}^{\overline{x}})^{q_2(i)}}{\lambda(a_{(i_0)}^{\overline{x}})^{|q_2|}}.$$ Note that $\lambda(a_{(i_0)}^{\overline{x}})=a_{i_0,i_0}^{\overline{x}}$ in the ${\cal O}_{\mathbb{C}_p}[[U]]$-algebra $\prod_J{\cal O}_{\mathbb{C}_p}[[U]]$ since ${\cal M}$ is standard normal. Let $\widetilde{\cal B}^r_{-}({\cal M}_{\overline{x}})$ be the matrix defined by the same recipe, but now using $\eta^{-1}$ in place of $\eta$.\\

\addtocounter{satz}{1}{\bf \arabic{section}.\arabic{satz}}\newcounter{specde1}\newcounter{specde2}\setcounter{specde1}{\value{section}}\setcounter{specde2}{\value{satz}} Let $\xi\in A$ be a unit, let $(s,t)\in{\cal F}$, let $r_1, r_2\in\mathbb{N}$, for $\ell=1$ and $\ell=2$ let ${I}^{(\ell)}$ be a countable index set, $i_0^{(\ell)}\in{I}^{(\ell)}$ an element and ${\cal M}_{\ell}$ a standard $1$-normal (with respect to $i_0^{(\ell)}$) nuclear ${I}^{(\ell)}\times{I}^{(\ell)}$-matrix over $\widehat{A}$. Arguing as in \ref{nucov}, where we proved that the matrices ${\cal B}^{r}({\cal M})$ are nuclear, we see that the trace$$g_{\overline{x},\xi,{\cal M}_1,{\cal M}_2}^{r_1,r_2,s,t}(U):=\tr_{}(\xi_{\overline{x}}^t\widetilde{\cal B}^{s-r_1}({\cal M}_{1,\overline{x}})\otimes\bigwedge^{r_1}({\cal M}_1)_{\overline{x}}\otimes\widetilde{\cal B}^{-s-r_2}_{-}({\cal M}_{2,\overline{x}})\otimes\bigwedge^{r_2}({\cal M}_2)_{\overline{x}})$$(the fibre $\xi_{\overline{x}}\in R$ is defined as in \arabic{lfnd1}.\arabic{lfnd2} by viewing $\xi$ as a $1\times1$-matrix) is a well defined element in ${\cal O}_{\mathbb{C}_p}[[U]]$, i.e. the infinite sum of diagonal elements of this tensor product matrix converges in ${\cal O}_{\mathbb{C}_p}[[U]]$. We may view it as a function on $D^{>0}$. Let $\nu\in\mathbb{Q}$ satisfy both \ref{lubtat} and the condition from \arabic{lide1}.\arabic{lide2} for both ${\cal M}_1$ and ${\cal M}_2$ so that we may form the matrices ${\cal B}^{s-r_1}({\cal M}_1)$ and ${\cal B}^{-s-r_2}_{-}({\cal M}_2)$ over $A\widehat{\otimes}_RB$ with $B=(B(\nu)_K)^0$. Recall the embedding $\iota:D^{\ge\nu}\to D^{>0}$ from \ref{lubtat} and that we view the free variable $V$ as standard coordinate on the source $D^{\ge\nu}$, and the free variable $U$ as standard coordinate on the target $D^{>0}$ of $\iota$. For a matrix ${\cal N}$ with coefficients in $A\widehat{\otimes}_RB$ and for $y\in D^{\ge\nu}$ we denote by ${\cal N}|_{V=y}$ the matrix with entries in $\widehat{A}$ obtained from ${\cal N}$ by specializing elements $a\otimes V^n\in A\widehat{\otimes}_RB$ (for $n\in\mathbb{N}_0$) to $a\otimes y^n\in \widehat{A}$.\\

\begin{lem}\label{spufas} For $K$-rational points $y\in D^{\ge \nu}$ we have$$g_{\overline{x},\xi,{\cal M}_1,{\cal M}_2}^{r_1,r_2,s,t}(\iota(y))=$$$$\tr(((\xi^t{\cal B}^{s-r_1}({\cal M}_1)\otimes\bigwedge^{r_1}({\cal M}_1)\otimes{\cal B}^{-s-r_2}_{-}({\cal M}_2)\otimes\bigwedge^{r_2}({\cal M}_2))|_{V=y})_{\overline{x}})$$\end{lem}

{\sc Proof:} Taking $\overline{x}$-fibres commutes with $\otimes$, thus$$\xi^t_{\overline{x}}({\cal B}^{s-r_1}({\cal M}_{1})|_{V=y})_{\overline{x}}\otimes\bigwedge^{r_1}({\cal M}_{1})_{\overline{x}}\otimes({\cal B}^{-s-r_2}_{-}({\cal M}_2)|_{V=y})_{\overline{x}}\otimes\bigwedge^{r_2}({\cal M}_2)_{\overline{x}}$$$$=((\xi^t{\cal B}^{s-r_1}({\cal M}_1)\otimes\bigwedge^{r_1}({\cal M}_1)\otimes{\cal B}^{-s-r_2}_{-}({\cal M}_2)\otimes\bigwedge^{r_2}({\cal M}_2))|_{V=y})_{\overline{x}}.$$Therefore it suffices to show$$({\cal B}^{r}({\cal M})|_{V=y})_{\overline{x}}=\widetilde{\cal B}^{r}({\cal M}_{\overline{x}})|_{U=\iota(y)}\quad\mbox{ and }\quad({\cal B}_{-}^{r}({\cal M})|_{V=y})_{\overline{x}}=\widetilde{\cal B}^{r}_{-}({\cal M}_{\overline{x}})|_{U=\iota(y)},$$ for standard $1$-normal nuclear matrices ${\cal M}$ over $\widehat{A}$ and $r\in\mathbb{Z}$. This is essentially a statement on commutation of the two operations "${\cal M}\mapsto{\cal B}^r({\cal M})$" and "taking the $f$-fold $\sigma$-power of a square matrix". In our situation this holds since ${\cal M}$ is standard normal, as we will now explain. For such ${\cal M}$ we keep the notation from \arabic{faslim1}.\arabic{faslim2}. From \ref{lubtat} it follows that $\widetilde{\cal B}^{r}({\cal M}_{\overline{x}})|_{U=\iota(y)}$ is the matrix constructed by the same recipe as $\widetilde{\cal B}^{r}({\cal M}_{\overline{x}})$, but using $$(a_{i_0,i_0}^{\overline{x}})^y=\exp(y\log(a_{i_0,i_0}^{\overline{x}}))\in R_f$$ in place of $\eta$. Observe that $a_{i_0,i_0}^{\overline{x}}=(a_{i_0,i_0})_{\overline{x}}$ where $(a_{i_0,i_0})_{\overline{x}}$ is defined as in \arabic{lfnd1}.\arabic{lfnd2} by viewing the $(i_0,i_0)$-entry $a_{i_0,i_0}$ of ${\cal M}$ as a $1\times 1$-matrix --- this is because ${\cal M}$ is standard. In particular we see $(a_{i_0,i_0}^{\overline{x}})^y=((a_{i_0,i_0})_{\overline{x}})^y\in R$. Let $(M,\phi)$ be the $\sigma$-module over $\widehat{A}$ such that the action of $\phi$ on a formal basis $\{e_i\}_{i\in I}$ of $M$ is given by ${\cal M}$. As in \arabic{modlim1}.\arabic{modlim2} consider the $\widehat{A}$-algebra $D=\sym_{\widehat{A}}(M)[\frac{1}{e_{i_0}}]^0$ of degree zero elements in $\sym_{\widehat{A}}(M)[\frac{1}{e_{i_0}}]$. Let $B^r(M)$ be its completion as in \arabic{modlim1}.\arabic{modlim2}. Denote by $\psi$ the natural $\sigma$-linear ring endomorphism of $B^r(M)$ defined by $\phi$, as in \arabic{modlim1}.\arabic{modlim2}. Then ${\cal B}^{r}({\cal M})|_{V=y}$ is the matrix of the $\sigma$-linear endomorphism $\psi_{y+r}=B^r(\phi)|_{V=y}=((a_{i_0,i_0})_{\overline{x}})^{y+r}\psi$ (use that ${\cal M}$ is standard). Hence $({\cal B}^{r}({\cal M})|_{V=y})_{\overline{x}}$ is the matrix of the $R_f$-linear endomorphism $(\psi_{y+r}^f)_{\overline{x}}$ which the $f$-fold iterate $\psi_{y+r}^f$ of $\psi_{y+r}$ induces on the fibre $B^r(M)_{\overline{x}}=B^r(M)\otimes_{\widehat{A}}R_f$ (formed with respect to the Teichm\"uller lift $x:\widehat{A}\to R_f$ of $\overline{x}$). On the other hand we can view $B^r(M)_{\overline{x}}$ as the completion (analogously to \arabic{modlim1}.\arabic{modlim2}) of $\sym_{R_f}(M_{\overline{x}})[\frac{1}{e_{i_0}}]^0$ (with $M_{\overline{x}}=M\otimes_{\widehat{A}}R_f$). Then $\widetilde{\cal B}^{r}({\cal M}_{\overline{x}})|_{U=\iota(y)}$ is the matrix of the $R_f$-linear endomorphism $((a_{i_0,i_0})_{\overline{x}})^{y+r}\psi_{f,\overline{x}}$ of $B^r(M)_{\overline{x}}$ where $\psi_{f,\overline{x}}$ is the $R_f$-linear ring endomorphism of $B^r(M)_{\overline{x}}$ induced by the endomorphism which the $f$-fold iterate $\phi^f$ of $\phi$ induces on $M_{\overline{x}}$. 
Thus it remains to show $(\psi_{y+r}^f)_{\overline{x}}=((a_{i_0,i_0})_{\overline{x}})^{y+r}\psi_{f,\overline{x}}$. Now we clearly have $(\psi_{y+r}^f)_{\overline{x}}=((a_{i_0,i_0})_{\overline{x}})^{y+r}(\psi^f)_{\overline{x}}$ with $(\psi^f)_{\overline{x}}$ the fibre of $\psi^f$ in $B^r(M)_{\overline{x}}$. Therefore we conclude using the functoriality $(\psi^f)_{\overline{x}}=\psi_{f,\overline{x}}$ of the ($\sigma$-linear) functor $\sym_{\widehat{A}}(?)$.\\

\addtocounter{satz}{1}{\bf \arabic{section}.\arabic{satz}} For $f\in\mathbb{N}$ let $T_f$ be the set of all closed points of $X$ of degree $f$. Let $A_f=A\otimes_RR_f$. Note that the $f$-fold $\sigma$-power $({\cal D}^{\wedge i})^{(\sigma)^f}$ (as defined in \arabic{lfnd1}.\arabic{lfnd2}) is the matrix describing the endomorphism which the $R_f$-algebra endomorphism $\sigma^f\otimes1$ of $A_f$ induces on $\Omega^i_{A_f/R_f}=\Omega^i_{A/R}\otimes_RR_f$. Therefore we may apply \ref{tracf} to the situation obtained by base change $\otimes_RR_f$, with $\sigma^f\otimes1\in{\rm End}(A_f)$ replacing $\sigma\in{\rm End}(A)$. We get that$$S_{\overline{x}}:=\sum_{0\le i\le d}(-1)^i\tr_{}(({\cal D}^{\wedge i})_{\overline{x}})$$for $\overline{x}\in T_f$ is invertible in $R_f$. For $0\le j\le d$ we may define$$h_{f,\xi,{\cal M}_1,{\cal M}_2}^{r_1,r_2,j,s,t}(U):=\sum_{\overline{x}\in T_f}\frac{\tr_{}(({\cal D}^{\wedge d-j})_{\overline{x}})}{S_{\overline{x}}}g^{r_1,r_2,s,t}_{\overline{x},\xi,{\cal M}_1,{\cal M}_2}\in{\cal O}_{\mathbb{C}_p}[[U]],$$$$D_{\xi,{\cal M}_1,{\cal M}_2}^{r_1,r_2,j,s,t}(T,U):=\exp(-\sum_{f=1}^{\infty}\frac{h_{f,\xi,{\cal M}_1,{\cal M}_2}^{r_1,r_2,j,s,t}(U)}{f}T^f)\in\mathbb{C}_p[[T,U]].$$

\begin{satz}\label{ganzfak} If ${\cal M}_1$ and ${\cal M}_2$ are $\sigma$-similar to $1$-normal nuclear overconvergent matrices over $A$, then $D_{\xi,{\cal M}_1,{\cal M}_2}^{r_1,r_2,j,s,t}(T,U)$ defines a holomorphic function on $\mathbb{A}_{\mathbb{C}_p}^1\times D^{>0}$. There exists a nuclear overconvergent matrix ${\cal N}$ over $A\widehat{\otimes}_RB$ which is $\sigma$-similar to$$\xi^t{\cal B}^{s-r_1}({\cal M}_1)\otimes\bigwedge^{r_1}({\cal M}_1)\otimes{\cal B}^{-s-r_2}_{-}({\cal M}_2)\otimes\bigwedge^{r_2}({\cal M}_2),$$and for $K$-rational points $y\in D^{\ge\nu}$ we have$$D_{\xi,{\cal M}_1,{\cal M}_2}^{r_1,r_2,j,s,t}(T,\iota(y))=\det(1-\psi[{\cal N}|_{V=y}\otimes{\cal D}^{\wedge d-j}]T).$$\end{satz}

{\sc Proof:} The existence of ${\cal N}$ follows from \ref{nucov} and \ref{funclim}. Next let us make a general remark. For a nuclear overconvergent matrix ${\cal M}$ over $A$ we defined the completely continuous operator $\psi[{\cal M}]=\psi_A[{\cal M}]$ in \arabic{dwrat1}.\arabic{dwrat2} relative to the Frobenius endomorphism $\sigma$ on $A$. Now consider the $f$-fold $\sigma$-power ${\cal M}^{(\sigma)^f}$ of ${\cal M}$ from \arabic{lfnd1}.\arabic{lfnd2} and {\it view it as a matrix over} $A_f=A\otimes_RR_f$. As such we define the $K_f=R_f\otimes\mathbb{Q}$-linear completely continuous operator $\psi_{A_f}[{\cal M}^{(\sigma)^f}]$ relative to the Frobenius endomorphism $\sigma^f$ on $A_f$. One finds $$\psi_{A_f}[{\cal M}^{(\sigma)^f}]=\psi_{A}[{\cal M}]^{f}\otimes_KK_f.$$We apply this to ${\cal M}={\cal N}|_{V=y}\otimes {\cal D}^{\wedge d-j}$ and obtain \begin{align}
\tr_K(\psi_A[{\cal N}|_{V=y}\otimes{\cal D}^{\wedge d-j}]^f) & = \tr_{K_f}(\psi_A[{\cal N}|_{V=y}\otimes{\cal D}^{\wedge d-j}]^f\otimes_KK_f) \notag \\
{} & =\tr_{K_f}(\psi_{A^f}[({\cal N}|_{V=y}\otimes{\cal D}^{\wedge d-j})^{(\sigma)^f}]) \notag \\
{} &  =\tr_{K_f}(\psi_{A^f}[({\cal N}|_{V=y})^{(\sigma)^f}\otimes({\cal D}^{\wedge d-j})^{(\sigma)^f}]) \notag \\
{} & =\sum_{\overline{x}\in T_f}\frac{\tr(({\cal D}^{\wedge d-j})_{\overline{x}})\tr(({\cal N}|_{V=y})_{\overline{x}})}{S_{\overline{x}}}. \notag
\end{align}where for the last equality we applied \ref{tracf}. But $$\tr(({\cal N}|_{V=y})_{\overline{x}})=\tr(((\xi^t{\cal B}^{s-r_1}({\cal M}_1)\otimes\bigwedge^{r_1}({\cal M}_1)\otimes{\cal B}^{-s-r_2}_{-}({\cal M}_2)\otimes\bigwedge^{r_2}({\cal M}_2))|_{V=y})_{\overline{x}})$$which by \ref{spufas} is equal to $g_{\overline{x},\xi,{\cal M}_1,{\cal M}_2}^{r_1,r_2,s,t}(\iota(y))$. Thus the stated formula is proven since its right hand side may be written as $$\exp(-\sum_{f=1}^{\infty}\frac{\tr_K(\psi[{\cal N}|_{V=y}\otimes{\cal D}^{\wedge d-j}]^f)}{f}T^f).$$Furthermore the points $\iota(y)$ for $K$-rational points $y\in D^{\ge\nu}$ are Zariski dense in $\iota(D^{\ge\nu})$, therefore we get the equality of holomorphic functions$$D_{\xi,{\cal M}_1,{\cal M}_2}^{r_1,r_2,j,s,t}(T,U)=\det(1-\psi[{\cal N}\otimes{\cal D}^{\wedge d-j}]T)$$on $D^{>0}\times\iota(D^{\ge\nu})$, where in the function on the right hand side we substitute $V$ by $\iota^{-1}(U)$. But the right hand side extends to a holomorphic function on $\mathbb{A}^1_{\mathbb{C}_p}\times\iota(D^{\ge\nu})$, since $\psi[{\cal N}\otimes{\cal D}^{\wedge d-j}]$ is completely continuous by \ref{nucl}. The definition of $D_{\xi,{\cal M}_1,{\cal M}_2}^{r_1,r_2,j,s,t}(T,U)$ and \ref{convext} now show the holomorphy on all of $\mathbb{A}_{\mathbb{C}_p}^1\times D^{>0}$, completing the proof.\\  

\addtocounter{satz}{1}{\bf \arabic{section}.\arabic{satz}} Let $\alpha\in\widehat{A}$ be a unit. For closed points $\overline{x}\in X$ define $\alpha_{\overline{x}}\in R$ as in \arabic{lfnd1}.\arabic{lfnd2} by viewing $\alpha$ as a $1\times 1$-matrix. For $\kappa\in\ho_{K\mbox{-}an}(R^{\times},\mathbb{C}_p^{\times})$ we ask for the twisted $L$-function$$L(\alpha,T,\kappa):=\prod_{\overline{x}\in X}\frac{1}{1-\kappa(\alpha_{\overline{x}})T^{\deg(\overline{x})}}.$$It can be written as a power series with coefficients in ${\cal O}_{\mathbb{C}_p}$, hence is trivially holomorphic on $D^{>0}$ (in the variable $T$).\\

\addtocounter{satz}{1}{\bf \arabic{section}.\arabic{satz}}
\newcounter{ordgeo1}\newcounter{ordgeo2}\setcounter{ordgeo1}{\value{section}}\setcounter{ordgeo2}{\value{satz}} We say that $\alpha\in\widehat{A}$ is {\it ordinary geometric} if there exists a nuclear $G\times G$-matrix ${\cal H}=(h_{g_1,g_2})_{g_1,g_2\in G}$ over $\widehat{A}$, a non negative integer $j\in\mathbb{N}_0$ and a nested sequence of $(j+1)$ finite subsets $G_0\subset G_1\subset\ldots\subset G_j$ of the (countable) index set $G$ such that:\\(i) ${\cal H}$ is $\sigma$-similar to a nuclear overconvergent matrix over $A$.\\(ii) $h_{g_1,g_2}=0$ whenever there is a $0\le\ell\le j$ with $g_2\in G_\ell$ and $g_1\notin G_\ell$. Thus, ${\cal H}$ is in block triangular form.\\(iii) $\pi^{\ell+1}$ divides $h_{g_1,g_2}$ whenever $g_2\notin G_\ell$, for all $0\le\ell\le j$.\\(iv) For all $0\le\ell\le j$ the element$$H_{\ell}:=\pi^{-\sum_{i=1}^{\ell}i(c_i-c_{i-1})}\det((h_{g_1,g_2})_{g_1,g_2\in G_{\ell}})$$of $\widehat{A}$ is a unit, where we set $c_{\ell}=|G_{\ell}|$. Set $H_{-1}=1$.\\(v) We have $\alpha=H_j/H_{j-1}=\pi^{-j(c_j-c_{j-1})}\det((h_{g_1,g_2})_{g_1,g_2\in (G_{j}-G_{j-1})})$.\\The meaning of this definition is that $\alpha$ is the determinant of the pure slope $j$ part (as a unit root $\sigma$-module) of a nuclear $\sigma$-module over $A$ which is ordinary up to slope $j$ and overconvergent (but neither the pure slope $j$ part itself nor its determinant need to be overconvergent). See \cite{wahrk} for details on the Hodge-Newton decomposition by slopes.\\ 

\begin{satz}\label{merwei} Suppose $\alpha$ is ordinary geometric. Then there exists a meromorphic function $L_\alpha$ on the $\mathbb{C}_p$-rigid space $\mathbb{A}_{\mathbb{C}_p}^1\times{\cal W}$ whose pullback to $\mathbb{A}_{\mathbb{C}_p}^1$ via $$\mathbb{A}_{\mathbb{C}_p}^1\to\mathbb{A}_{\mathbb{C}_p}^1\times{\cal W},\quad t\mapsto(t,\kappa),$$for any $\kappa\in\ho_{K\mbox{-}an}(R^{\times},\mathbb{C}_p^{\times})={\cal W}(\mathbb{C}_p)$ is a continuation of $L(\alpha,T,\kappa)$.\end{satz}

{\sc Proof:} We treat every component ${\cal W}_{(s,t)}$ of ${\cal W}$ separately, so let us fix $(s,t)\in{\cal F}$. Keeping the notation from \arabic{ordgeo1}.\arabic{ordgeo2} we begin with some definitions. For $0\le\ell\le j$ let $I^{({\ell})}$ be the index set of the nuclear matrix $\bigwedge^{c_{\ell}}({\cal H})$. Our assumptions on ${\cal H}$ imply that $\pi^{-j(c_j-c_{j-1})}\bigwedge^{c_{\ell}}({\cal H})$ is standard normal with respect to some $i^{({\ell})}_0\in I^{({\ell})}$. Moreover it is $\sigma$-similar to a nuclear overconvergent $I^{({\ell})}\times I^{({\ell})}$ matrix. Thus we may apply \ref{discrdec} to get a $\xi_{\ell}\in A$ and a $1$-normal (with respect to $i^{({\ell})}_0$) nuclear overconvergent $I^{({\ell})}\times I^{({\ell})}$-matrix ${\cal M}_{\ell}$ over $A$ such that $\xi_{\ell}^{q-1}$ is $\sigma$-similar to $1\in A$, and $\xi_{\ell}{\cal M}_{\ell}$ is $\sigma$-similar to $\pi^{-j(c_j-c_{j-1})}\bigwedge^{c_{\ell}}({\cal H})$. By \ref{hndec} there is a standard $1$-normal (with respect to $i^{({\ell})}_0$) nuclear $I^{({\ell})}\times I^{({\ell})}$-matrix ${\cal M}'_{\ell}$ over $\widehat{A}$ which is $\sigma$-similar to ${\cal M}_{\ell}$. Let ${\cal M}'_{\ell,unit}\in\widehat{A}$ be the $(i^{({\ell})}_0,i^{({\ell})}_0)$-entry of ${\cal M}'_{\ell}$. This is a $1$-unit. Then $\xi_{\ell}{\cal M}'_{\ell,unit}\in\widehat{A}$ is $\sigma$-similar to the $(i^{({\ell})}_0,i^{({\ell})}_0)$-entry of $\pi^{-j(c_j-c_{j-1})}\bigwedge^{c_{\ell}}({\cal H})$ which we denote by $a_{\ell}$. We will need these definitions for $\ell=j-1$ and $\ell=j$ if $j>0$. If $j=0$ we set $\xi_{-1}={\cal M}_{-1}={\cal M}'_{-1}={\cal M}'_{-1,unit}=a_{-1}=1\in R$. Our definitions imply $\alpha=a_{j}/a_{j-1}$, thus if we set $\xi=\xi_j/\xi_{j-1}$ and $\mu={\cal M}'_{j,unit}/{\cal M}'_{j-1,unit}$ we find that $\alpha$ is $\sigma$-similar to $\xi\mu$. Let$$H(T,U)=\prod_{\overline{x}\in X}\frac{1}{1-\xi_{\overline{x}}^t\mu_{\overline{x}}^s(1+F([\pi^m\log(\mu_{\overline{x}})]))T^{\deg{\overline{x}}}}.$$This is a holomorphic function on $D^{>0}\times D^{>0}$ where we view $T$ (resp. $U$) as coordinate for the first (resp. second) factor $D^{>0}$. Recall from \arabic{weicons1}.\arabic{weicons2} the finite \'{e}tale covering of rigid spaces $D^{>0}\to{\cal W}_{(s,t)}$ which on $\mathbb{C}_p$-valued points is given by$$D^{>0}\to{\cal W}_{(s,t)}(\mathbb{C}_p)\cong\ho_{K\mbox{-}an}(\overline{U}_{R}^{(1)},\mathbb{C}_p^{\times})$$$$z\longmapsto[\overline{w}\mapsto1+F([\pi^m\log(\overline{w})])(z)].$$We see that for any $w\in{U}_{R}^{(1)}$ the holomorphic function $1+F([\pi^m\log({w})])(U)$ in the variable $U$ on $D^{>0}$ descends to ${\cal W}_{(s,t)}$. Thus our $H(T,U)$ descends to a holomorphic function $L_{\alpha,(s,t)}$ on $D^{>0}\times{\cal W}_{(s,t)}$. Moreover, for $\kappa\in{\cal W}_{(s,t)}(\mathbb{C}_p)\subset{\cal W}(\mathbb{C}_p)=\ho_{K\mbox{-}an}({R}^{\times},\mathbb{C}_p^{\times})$ the pullback of $L_{\alpha,(s,t)}$ via$$D^{>0}\to D^{>0}\times{\cal W}_{(s,t)},\quad t\mapsto(t,\kappa)$$is $L(\alpha,T,\kappa)$: this is immediate since $\xi_{\overline{x}}\mu_{\overline{x}}=\alpha_{\overline{x}}$ is the decomposition of $\alpha_{\overline{x}}\in R^{\times}$ according to $R^{\times}=\mathbb{\mu}_{q-1}\times U^{(1)}_R$, for any $\overline{x}\in X$. These considerations also show that for $K$-rational points $y\in D^{\ge\nu}$ we have \begin{gather}H(T,\iota(y))=L(\xi^t\mu^s\mu^y,T)\tag{$1$}\end{gather}with $\iota: D^{\ge\nu}\to D^{>0}$ from \ref{lubtat}. To show that $L_{\alpha,(s,t)}$ is meromorphic on $\mathbb{A}_{\mathbb{C}_p}^1\times{\cal W}_{(s,t)}$ it is enough to show that $H(T,U)$ is meromorphic on $\mathbb{A}_{\mathbb{C}_p}^1\times D^{>0}$. Consider the $\mathbb{C}_p[[T,U]]$-element$${\underline H}(T,U)=\prod_{r_1,r_2\ge1}(\prod_{i=0}^{d}D_{\xi,{\cal M}'_j,{\cal M}'_{j-1}}^{r_1,r_2,d-i,s,t}(T,U)^{(-1)^{i-1}})^{(-1)^{r_1+r_2}r_1r_2}.$$ By \ref{ganzfak} each factor $D_{\xi,{\cal M}'_j,{\cal M}'_{j-1}}^{r_1,r_2,d-i,s,t}(T,U)$ is holomorphic on $\mathbb{A}_{\mathbb{C}_p}^1\times D^{>0}$. Moreover, since $\bigwedge^{r_{\ell}}({\cal M}'_{\ell})$ is divisible by $\pi^{r_{\ell}-1}$ it also follows from \ref{ganzfak} that $\ord_{\pi}(1-D_{\xi,{\cal M}'_j,{\cal M}'_{j-1}}^{r_1,r_2,d-i,s,t}(T,U))$ tends to infinity as $r_1+r_2$ tends to infinity (if the index set $G$ is finite then the above product is even finite). Therefore ${\underline H}(T,U)$ is meromorphic on $\mathbb{A}_{\mathbb{C}_p}^1\times D^{>0}$. We claim $${\underline H}(T,U)=H(T,U)$$as meromorphic functions on $D^{>0}\times D^{>0}$. As in \ref{ganzfak} it is enough to check this on all subsets $D^{>0}\times{\iota(y)}\subset D^{>0}\times D^{>0}$ for $K$-rational points $y\in D^{\ge\nu}$. From \ref{rk1res} we get the following equalities in the Grothendieck group $\Delta(\widehat{A})$:$$[\xi^t({\cal M}'_{j,unit})^s({\cal M}'_{j,unit})^y]=\bigoplus_{r\ge1}(-1)^{r-1}r[\xi^t{\cal B}^{s-r}({\cal M}'_j)|_{V=y}\otimes\bigwedge^r({\cal M}'_j)]$$$$[({\cal M}'_{j-1,unit})^{-s}({\cal M}'_{j-1,unit})^{-y}]=\bigoplus_{r\ge1}(-1)^{r-1}r[{\cal B}^{-s-r}_{-}({\cal M}'_{j-1})|_{V=y}\otimes\bigwedge^r({\cal M}'_{j-1})]$$(with the notation $|_{V=y}$ explained in \arabic{specde1}.\arabic{specde2} still in force: $V$ is the standard coordinate on $D^{\ge\nu}$). Together$$[\xi^t\mu^s\mu^y]=[\frac{\xi^t({\cal M}'_{j,unit})^s({\cal M}'_{j,unit})^y}{({\cal M}'_{j-1,unit})^s({\cal M}'_{j-1,unit})^y}]=$$$$\bigoplus_{r_1,r_2\ge1}(-1)^{r_1+r_2}r_1r_2[\xi^t({\cal B}^{s-r_1}({\cal M}'_j)\otimes\bigwedge^{r_1}{\cal M}'_j\otimes{\cal B}^{-s-r_2}_{-}({\cal M}'_{j-1})\otimes\bigwedge^{r_2}{\cal M}'_{j-1})|_{V=y}].$$Combining with \ref{ganzfak} and the trace formula \ref{tracf} we get \begin{gather}{\underline H}(T,\iota(y))=L(\xi^t\mu^s\mu^y,T).\tag{$2$}\end{gather}Comparing $(1)$ and $(2)$ completes the proof.\\

\section{Higher rank}

\addtocounter{satz}{1}{\bf \arabic{section}.\arabic{satz}} A finite rank $\sigma$-module $(M,\phi)$ over $\widehat{A}$ is called {\it ordinary} if it admits a separated and exhausting $\phi$-stable filtration by free sub-$\widehat{A}$-modules $$0=M_0\subset M_1\subset M_2\subset\ldots$$of $M$ with free quotients, such that each quotient $(M_i/M_{i+1},\phi)$ is of the form $(U_i,\pi^i.\phi_i)$ where $(U_i,\phi_i)$ is a unit root $\sigma$-module; that is, $\widehat{A}_{\sigma}\otimes U_i\to U_i, a\otimes u\mapsto a.\phi_i(u)$ is bijective. The $(U_i,\phi_i)$ are called the {\it graded pieces of} $(M,\phi)$, and $(U_0,\phi_0)$ is called the {\it unit root part of} $(M,\phi)$, also denoted by $\phi_{unit}$.\\

\begin{satz}\label{hndecov} (Hodge-Newton slope decomposition for overconvergent $\sigma$-modules) Let ${\cal M}$ be the matrix, in some basis, of a graded piece of an ordinary and overconvergent finite rank $\sigma$-module $(M,\phi)$ over $\widehat{A}$. Then there exists a convergent series representation$$[{\cal M}]=\sum_{r\ge 1}\pm[{\cal C}_r]$$in $\Delta(\widehat{A})$ with nuclear overconvergent matrices ${\cal C}_r$ over $A$.
\end{satz}

{\sc Proof:} (1) By induction on $m$ we prove that for each $m\in\mathbb{N}_0$ there exist finite index sets $J_m^1, J_m^2$, ordinary overconvergent $\sigma$-modules $\alpha_t$ and $\beta_t$ of finite rank for each $t\in J_m^1\cup J_m^2$, with $(\beta_t)_{unit}$ of rank one, and integers $m_t\ge m$ for each $t\in J_m^2$, such that\begin{gather}[{\cal M}]=(\sum_{t\in J_m^2}\pm[\pi^{m_t}(\alpha_t)_{unit}\otimes(\beta_t)^{-1}_{unit}])+(\sum_{t\in J_m^1}\pm[\alpha_t\otimes(\beta_t)^{-1}_{unit}])\tag{$*$}\end{gather}in $\Delta(\widehat{A})$. Here, by abuse of notation, we identify a finite rank $\sigma$-module with the $\sigma$-similarity class of matrices it corresponds to. For $m=0$ one has $[{\cal M}]=[\alpha_{unit}\otimes\beta_{unit}^{-1}]$ for some $\alpha, \beta$, by \cite{wahrk} 6.2. Now let us pass from $m$ to $m+1$. Fix $t\in J_m^2$. Let $(\alpha_{t,t'})_{t'\in T_t}$ be the set of  higher graded pieces of ${\alpha}_t$ (i.e. $(\alpha_t)_{unit}$ omitted). By \cite{wahrk} 6.2 again, there exist for each $t'\in T_t$ ordinary overconvergent finite rank $\sigma$-modules $\widetilde{\alpha}_{t,t'}$ and $\widetilde{\beta}_{t,t'}$, with $(\widetilde{\beta}_{t,t'})_{unit}$ of rank one, such that $[\alpha_{t,t'}]=[(\widetilde{\alpha}_{t,t'})_{unit}\otimes(\widetilde{\beta}_{t,t'})^{-1}_{unit}]$. Thus$$[({\alpha}_{t})_{unit}\otimes({\beta}_{t})^{-1}_{unit}]=[{\alpha}_{t}\otimes({\beta}_{t})^{-1}_{unit}]-\sum_{t\in T_t}[\pi^{m_{t'}}({\alpha}_{t,t'})\otimes(\beta_t)^{-1}_{unit}]$$$$=[{\alpha}_{t}\otimes({\beta}_{t})^{-1}_{unit}]-\sum_{t\in T_t}[\pi^{m_{t'}}(\widetilde{\alpha}_{t,t'})_{unit}\otimes(\beta_t\otimes\widetilde{\beta}_{t,t'})^{-1}_{unit}]$$with integers $m_{t'}\ge1$ (the higher slopes of ${\alpha_t}$). Inserting this into the formula given by induction hypothesis for $m$ gives the formula for $m+1$.\\(2) To get the desired convergent series representation for ${\cal M}$ it is now enough to express, for $t\in J_m^1$, the terms $({\beta}_{t})^{-1}_{unit}$ in $(*)$ through overconvergent matrices. This is achieved by factoring ${\beta}_{t}$ according to \ref{discrdec} and applying \ref{rk1res} (with $V=0$ and $s=-1$ there) to the $1$-normal overconvergent factor of ${\beta}_{t}$.\\

\addtocounter{satz}{1}{\bf \arabic{section}.\arabic{satz}} As a corollary of Theorem \ref{hndecov} (and \ref{lfndbkt}) we recover Wan's result: that $L({\cal M},T)$ is a meromorphic function on $\mathbb{A}^1$. 


\end{document}